\documentclass[preprint]{elsarticle} 

\makeatletter
\def\ps@pprintTitle{%
\let\@oddhead\@empty
\let\@evenhead\@empty
\def\@oddfoot{\centerline{\thepage}}%
\let\@evenfoot\@oddfoot}
\makeatother

\usepackage{etoolbox}
\patchcmd{\MaketitleBox}{\footnotesize\itshape\elsaddress\par\vskip36pt}{\footnotesize\itshape\elsaddress\par\parbox[b][36pt]{\linewidth}{\vfill\hfill\textnormal{\today}\hfill\null\vfill}}{}{}%
\patchcmd{\pprintMaketitle}{\footnotesize\itshape\elsaddress\par\vskip36pt}{\footnotesize\itshape\elsaddress\par\parbox[b][36pt]{\linewidth}{\vfill\hfill\textnormal{\today}\hfill\null\vfill}}{}{}%

\usepackage[
colorlinks=true,%
breaklinks=true,%
linkcolor=blue,
urlcolor=blue,%
citecolor=blue,%
pdftitle={Global sensitivity analysis with limited data via D-MORPH regression: Application to char combustion}, 
pdfkeywords={Global sensitivity analysis, ANOVA, polynomial dimensional decomposition, surrogate modeling, D-MORPH regression, char combustion},	
pdfauthor={Dongjin Lee, Elle Lavichant, Boris Kramer},		
bookmarksopen=false,
pdfpagemode=UseNone]{hyperref}
\usepackage{algorithm}
\usepackage{algorithmicx}
\usepackage{algpseudocode}
\usepackage{tabularx}
\usepackage{array}
\usepackage{multirow}
\usepackage{threeparttable}
\usepackage{booktabs}
\usepackage{graphicx}
\usepackage{arydshln}
\usepackage{setspace}   
\usepackage{xcolor}
\usepackage[mathlines,running]{lineno}
\usepackage{epstopdf}
\usepackage{natbib} 

\usepackage[margin = 1.2in]{geometry}
\usepackage[T1]{fontenc}
\usepackage[english]{babel}
\usepackage[latin1]{inputenc}
\usepackage{mathtools}
\usepackage{amsmath}
\usepackage{amsfonts}
\usepackage{amsthm}
\usepackage{amssymb}
\usepackage{array}
\usepackage{subcaption}
\usepackage{siunitx}

\usepackage{color}
\usepackage{url}
\usepackage{cleveref}
\usepackage{graphicx}
\usepackage{breakcites}
\usepackage{soul} 
\usepackage{enumitem}
\usepackage[title,titletoc,toc]{appendix}

\usepackage[textsize=tiny]{todonotes}

\newcommand{\RNum}[1]{\uppercase\expandafter{\romannumeral #1\relax}}
\usepackage{amsthm}
\makeatletter
\def\els@aparagraph[#1]#2{\elsparagraph[#1]{#2\@addpunct{.}}}
\def\els@bparagraph#1{\elsparagraph*{#1\@addpunct{.}}}
\makeatother

{\noindent {\textbf{Proof}:} }%
{\hfill $\Box$ \\[1ex] }

\newcommand{\bit}{\begin{itemize}}
	\newcommand{\eit}{\end{itemize}}
\newcommand{\ben}{\begin{enumerate}}
	\newcommand{\een}{\end{enumerate}}

\newcommand {\real} {\mathbb{R}}

\DeclareMathOperator*{\argmin}{arg\,min}%
\newcommand{\rd}{\text{\upshape d}} 

\newcommand{\bA}{\ensuremath{\mathbf{A}}}

\newcommand{\bD}{\ensuremath{\mathbf{D}}}
\newcommand{\bE}{\ensuremath{\mathbf{E}}}
\newcommand{\bF}{\ensuremath{\mathbf{F}}}

\newcommand{\bH}{\ensuremath{\mathbf{H}}}
\newcommand{\bI}{\ensuremath{\mathbf{I}}}

\newcommand{\bK}{\ensuremath{\mathbf{K}}}

\newcommand{\bM}{\ensuremath{\mathbf{M}}}

\newcommand{\bR}{\ensuremath{\mathbf{R}}}

\newcommand{\bT}{\ensuremath{\mathbf{T}}}

\newcommand{\bW}{\ensuremath{\mathbf{W}}}
\newcommand{\bX}{\ensuremath{\mathbf{X}}}

\newcommand{\ba}{\ensuremath{\mathbf{a}}}
\newcommand{\bb}{\ensuremath{\mathbf{b}}}
\newcommand{\bc}{\ensuremath{\mathbf{c}}}

\newcommand{\bj}{\ensuremath{\mathbf{j}}}

\newcommand{\bu}{\ensuremath{\mathbf{u}}}
\newcommand{\bv}{\ensuremath{\mathbf{v}}}

\newcommand{\bx}{\ensuremath{\mathbf{x}}}

\newcommand {\bPhi} {\mbox{\boldmath $\Phi$}}

\newcommand{\cB}{\ensuremath{\mathcal{B}}}

\newcommand{\cF}{\ensuremath{\mathcal{F}}}

\newcommand{\cM}{\ensuremath{\mathcal{M}}}

\graphicspath{{./figures/}}
\begin{document}
\begin{frontmatter}
		
\title{{Global sensitivity analysis with limited data via sparsity-promoting D-MORPH regression: Application to char combustion}}
		
		
\author[hyu]{Dongjin Lee\corref{cor1}}
\ead{dlee46@hanyang.ac.kr}
\author[ucsd]{Elle Lavichant}
\author[ucsd]{Boris Kramer}
		
\cortext[cor1]{Corresponding author}
\address[hyu]{{\color{black}{Department of Automotive Engineering, Hanyang University, Seoul, South Korea}}}		
\address[ucsd]{Department of Mechanical and Aerospace Engineering, University of California San Diego, CA, United States}

\begin{abstract}
In uncertainty quantification, variance-based global sensitivity analysis quantitatively determines the effect of each input random variable on the output by partitioning the total output variance into contributions from each input. However, computing conditional expectations can be prohibitively costly when working with expensive-to-evaluate models. 
Surrogate models can accelerate this, yet their accuracy depends on the quality and quantity of training data, which is expensive to generate (experimentally or computationally) for complex engineering systems. Thus, methods that work with limited data are desirable. 
We propose a diffeomorphic modulation under observable response preserving homotopy (D-MORPH) regression to train a polynomial dimensional decomposition surrogate of the output that minimizes the number of training data. The new method first computes a sparse Lasso solution and uses it to define the cost function. A subsequent D-MORPH regression minimizes the difference between the D-MORPH and Lasso solution. The resulting D-MORPH \textcolor{black}{based} surrogate is more robust to input variations and more accurate with limited training data. 
We illustrate the accuracy and computational efficiency of the new surrogate for global sensitivity analysis using mathematical functions and an expensive-to-simulate model of char combustion. The new method is highly efficient, requiring only 15\% of the training data compared to conventional regression.  
\end{abstract}

		
\end{frontmatter}
\section{Introduction} \label{sec:intro}
Global sensitivity analysis is a powerful tool in uncertainty quantification to measure the influence of input variables on the quantity of interest (QoI), or output, of a simulation model. Global sensitivity analysis helped scientists and engineers identify the most influential input variables and to make better decisions concerning system design and operation in a broad range of applications, including space weather~\citep{jivani2023global,issan2023bayesian}, biological systems~\citep{kiparissides2009global,linden2022bayesian}, structural engineering~\citep{kala2016global,rahman2016surrogate}, or wind energy~\citep{carta2020global}. In contrast to local sensitivity analysis that determines the gradients of the output with respect to variations in each input, the global sensitivity analysis explores the full range of the input variables, providing a more complete picture~\citep{opgenoord2016,rahman2011global}.
Consider, for instance, an engineering system (e.g., a power plant) with three operating or design inputs that are mapped via a simulation to a QoI. When implementing quantitative global sensitivity analysis, one may find that one input possesses a sensitivity of 90\% to the QoI, whereas the combined sensitivity index of the other two inputs contributes only 10\% to the variation of the QoI. 
Methods to optimize the system or quantify its uncertainty, often scale badly with the number of input variables due to the curse of dimensionality~\citep{lee2023high, lee2023bi}. Global sensitivity analysis can assist the system designer by providing quantitative reasons to neglect the two parameters that only contribute 10\% to variations in the QoI and focusing only on the most influential parameter.

Variance-based sensitivity analysis (also called ANalysis Of VAriance (ANOVA) or Sobol sensitivity analysis~\citep{sobol2001global}) estimates the sensitivity of the QoI to each input variable by partitioning the total variance of the output into contributions from each input variable~\citep{opgenoord2016}. 
Variance-based methods are formulated as conditional variances and can be evaluated by Monte Carlo simulation or Latin hypercube sampling~\citep{helton2003latin}, see also the review~\citep{sudret2008global}. These sampling methods often face difficulties when applied to a computationally intensive model. Numerous global sensitivity analysis studies therefore leverage efficient surrogate models. 
The authors in~\citep{qian2018multifidelity} present a multi-fidelity framework for global sensitivity analysis that combines a high-fidelity model with multiple reduced-order models to efficiently provide unbiased estimates of global sensitivity measures. 
Commonly used surrogates are polynomial chaos expansion (PCE)~\citep{sudret2008global}, polynomial dimensional decomposition (PDD)~\citep{rahman2011global}, Gaussian process or Kriging~\citep{wang2013application}, support vector machine~\citep{cheng2017global}, and artificial neural network~\citep{li2016accelerate}. \textcolor{black}{PCE, in particular, has been efficiently used for global sensitivity analysis by integrating robust regression methods like weighted $l^1$-minimization~\citep{peng2014weighted}, reweighted $l^1$-minimization~\citep{yang2013reweighted}, and similar approaches~\citep{doostan2009least, tsilifis2019compressive}. However, these $l^1$-minimization techniques have their limitations in that the maximum number of non-zero terms cannot exceed the number of available data points. This becomes a significant constraint when dealing with limited data, especially since the number of PCE terms increases exponentially.} 
PDD, which is a Fourier-polynomial expansion of lower-variate component functions~\citep{rahman2008}, is very effective for global sensitivity analysis. The PDD surrogate truncates its bases effectively in a dimension-wise manner which, to some effect, alleviates the curse of dimensionality compared to PCE. The PDD surrogate, similar to ANOVA, partitions the variance of a QoI among its different input variables, enabling variance-based sensitivity analysis~\citep{rahman2011global}. 
The accuracy of PDD depends on the number of training samples. When using conventional regression methods (e.g., standard least squares) to compute the surrogate, the number of training samples must exceed the number of to-be-learned expansion coefficients. Otherwise, we need to reduce the number of PDD bases, which degrades its accuracy relative to the QoI. As is often the case for nonlinear QoIs, the PDD surrogate can require several hundreds to thousands of training samples. Given that each training sample is obtained via an expensive simulation (e.g., in our application of char combustion in Section~\ref{sec:5}, one simulation takes 7.9 CPU hours) training the PDD model becomes impractical.
Instead of reducing the basis size, one can solve an underdetermined linear system to train the surrogate model. The authors in~\citep{li2010,li2012} introduce the diffeomorphic modulation under observable response preserving homotopy (D-MORPH) regression. 
The D-MORPH regression finds a solution to the linear system, and hence a surrogate model, that exactly fits the training data and minimizes the variance of the surrogate predictions~\citep{li2010}. 
When minimizing the cost function, D-MORPH assigns weights to the unknown regression solution. The resulting D-MORPH-based solution is highly sensitive to this weight selection.  
Although the authors in~\citep{li2010,li2012} provide strategies for selecting these weights, their methods are specific to a particular problem and require generalization to be applicable across different problems. 
In \citep{lee2020practical}, the authors use a D-MORPH regression to train the generalized polynomial chaos expansion surrogate for dependent inputs. The D-MORPH method partitions the surrogate basis functions into two groups. The first group consists of fewer lower-order basis functions than the number of training samples, which ensures the accuracy of the solution. The second group consists of the higher-order basis functions, which are assumed to have expansion coefficients with smaller magnitude compared to those in the primary group. However, this assumption may not hold in practical problems, requiring a more general version of the D-MORPH regression.
In \citep{he2022adaptive} a least absolute shrinkage and selection operator (Lasso) regression is used to train a PDD surrogate. The Lasso regression includes a regularization term that penalizes the $l^1$ norm of the regression coefficients. While the Lasso regression effectively yields sparse solutions for underdetermined systems, it induces a bias into the QoI estimates in an effort to decrease the variance. Moreover, the Lasso regression can only infer as many non-zero coefficients as there are training samples. 

In this work, we develop a new regression method to compute a PDD surrogate in the limited data setting, i.e., where the linear systems for the surrogate coefficients are underdetermined. 
In the new method, we first compute a sparse---yet biased---Lasso solution to the underdetermined system for the PDD coefficients. We then define a new cost function that represents the difference between the D-MORPH and the Lasso-based estimates. \textcolor{black}{The D-MORPH framework iteratively minimizes the cost function with an additional term that represents the difference between a potential solution derived from the current iteration of D-MORPH and the previously obtained D-MORPH solution. Although the cost function seeks to find a sparse solution, this iterative process can yield values near zero, which could be found in the true solution. }
In leveraging both the D-MORPH and Lasso regression, we obtain a surrogate that is more robust to variations in inputs and produces more accurate solutions with limited training samples.   
%

We demonstrate the accuracy and the computational efficiency of the new regression method for global sensitivity analysis using mathematical functions. We evaluate the new regression method for an expensive-to-simulate char combustion simulation model. This simulation is based on the Eulerian framework for gas behavior and the Lagrangian framework for solid particle behavior and includes  chemical kinetics models.
A single simulation requires 7.9 CPU hours in a parallelized implementation with 15 CPUs. 
Through the new method, we reduce the necessary training samples to only 15\% of the data required for a conventional regression. 

The paper is organized as follows. Section~\ref{sec:2} covers the theoretical background, the problem setting, and briefly introduces global sensitivity analysis, the PDD surrogate, and the original D-MORPH regression. Section~\ref{sec:3} proposes a novel D-MORPH regression to train the PDD surrogate. Therein, we demonstrate the proposed method using two mathematical functions. Section~\ref{sec:4} uses the proposed method for global sensitivity analysis of char combustion, where we only need 151 samples for training the surrogate. In Section~\ref{sec:5}, we draw the conclusions and point to future research.

\section{Theoretical background} \label{sec:2} 
We present the preliminaries and define the input random variables in Section~\ref{sec:2.2} and output random variables in Section~\ref{sec:2.3}. Section~\ref{sec:2.4} presents the global sensitivity analysis method. We  summarize the polynomial dimensional decomposition surrogate in Section~\ref{sec:2.5} and the original D-MORPH regression in Section~\ref{sec:2.6}.   
%
\subsection{Input random variables} \label{sec:2.2} 
Let $\real$, and $\real_{0}^{+}$ denote the real numbers and non-negative real numbers, respectively. For a positive integer $N$, denote by $\mathbb{A}^N \subseteq \real^N$ a bounded or unbounded sub-domain of $\real^N$. 
Let $(\Omega,\cF,\mathbb{P})$ be an abstract probability space, with a sample space $\Omega$, a $\sigma$-algebra $\cF$ on $\Omega$, and a probability measure $\mathbb{P}:\cF\to[0,1]$. Consider an $N$-dimensional random vector $\bX:=(X_{1},\ldots,X_{N})^\intercal:\Omega \rightarrow \mathbb{A}^N$ that models the input uncertainties. We refer to $\bX$ as the random input vector or the input random variables. 
Denote by $F_{\bX}({\bx}):=\mathbb{P}\big[\cap_{i=1}^{N}\{ X_i \le x_i \}\big]$ the joint distribution function of $\bX$. In this work, we assume that the input random variables in $\bX$ are independent so that the joint probability density function is $f_{\bX}({\bX}):=\prod_{k=1}^{k=N}f_k(X_k)$. Here, $f_k(X_k)$ is the marginal probability density function of $X_k$ defined on the probability space $(\Omega_k,\mathcal{F}_k,\mathbb{P}_k)$. For $(\Omega,\cF,\mathbb{P})$, the image probability space is $(\mathbb{A}^N,\cB^{N},f_{\bX}(\bx)\rd\bx)$, where $\mathbb{A}^N$ is the image of $\Omega$ under the mapping $\bX:\Omega \to \mathbb{A}^N$ and $\cB^N:=\cB(\mathbb{A}^N)$ is the Borel $\sigma$-algebra on $\mathbb{A}^N\subset \mathbb{R}^N$.

\subsection{Output random variables}  \label{sec:2.3}
Given an input random vector $\bX$
with a known probability density function $f_{\bX}({\bx})$ on $\mathbb{A}^N \subseteq \real^N$, denote by $Y=y(\bX)$ a real-valued, square-integrable  function on $(\Omega, \cF)$.   In this work, we assume that the output function (or quantity of interest) $y$ belongs to the weighted $L^2$ space $$\left\{y:\mathbb{A}^N \to \real:~
\int_{\mathbb{A}^N} \left| y(\bx)\right|^2 f_{\bX}({\bx})\rd\bx < \infty \right\}.$$ If there is more than one output, then each component is associated with a measurement function $y_i$, $i=1,2,\ldots$.  The generalization for a multivariate output random vector is straightforward.

\subsection{Global sensitivity analysis} \label{sec:2.4} 
We review variance-based methods for global sensitivity analysis followed by a brief explanation of how to use the polynomial dimensional decomposition surrogate for variance-based sensitivity analysis.

\subsubsection{ANalysis Of VAriance (ANOVA) dimensional decomposition}\label{sec:2.4:1}
The ANOVA dimensional decomposition assesses the effect of inputs $\bX$ on an output $y(\bX)$. 
The dimensional decomposition can be expressed as 
\begin{equation}
\begin{split}
y(\bX)=&y_0 + \sum_{i=1}^N y_i(X_i) + \sum_{i_1=1}^{N-1}\sum_{i_2=i_1+1}^N y_{i_1 i_2}(X_{i_1},X_{i_2})+\cdots\\
&+\sum_{i_1=1}^{N-s+1}\cdots\sum_{i_s=i_{s-1}+1}^N y_{i_1\cdots i_s}(X_{i_1},\ldots,X_{i_s})+\cdots+y_{12\cdot N}(X_1,\ldots,X_N),
\end{split} 
\label{ANOVA2}
\end{equation} 
where the component functions $y_{i_1,\ldots,i_s}(x_{i_1},\ldots,x_{i_s})$ for $1\leq i_1 < \cdots < i_s \leq N$ and $s=1,\ldots,N$ are defined as
\begin{align*}
y_0 &:=\int_{\mathbb{R}^N}y(\bx)f_{\bX}(\bx){\rm{d}}\bx, \\
y_i(x_i)&:=\int_{\mathbb{R}^{N-1}}y(\bx)\prod_{j\neq i} f_{X_j}(x_j){\rm{d}}x_j-y_0, \\
y_{i_1 i_2}(x_{i_1},x_{i_2}) &:=\int_{\mathbb{R}^{N-2}}y(\bx)\prod_{j\neq \{i_1,i_2\}} f_{X_j}(x_j){\rm{d}}x_j-y_{i_1}(x_{i_1})-y_{i_2}(x_{i_2})-y_0,\\
y_{i_1\cdots i_s}(X_{i_1},\ldots,X_{i_s})&:=\int_{\mathbb{R}^{N-s}}y(\bx)\prod_{j\neq \{i_1,\ldots,i_s\}} f_{X_j}(x_j){\rm{d}}x_j-
\sum_{j_1<\cdots <j_{s-1} \subset \{i_1,\ldots,i_s\}} y_{j_1\cdots j_{s-1}}(x_{j_1},\ldots,x_{j_s-1})\\
& \ \  -\sum_{j_1<\cdots <j_{s-2} \subset \{i_1,\ldots,i_s\}} y_{j_1\cdots j_{s-2}}(x_{j_1},\ldots,x_{j_s-2}) -\cdots-\sum_{j \subset \{i_1,\ldots,i_s\}} y_{j}(x_{j}) -y_0.
\end{align*}
These component functions satisfy orthogonal properties such that
\begin{equation} \label{ortho}
\int_{\mathbb{R}^N}y_{i_1\cdots i_s}(x_{i_1},\ldots,x_{i_s})f_{\bX}(\bx)\rm{d}\bx=0, \qquad
\int_{\mathbb{R}^N}y_{i_1\cdots i_s}(x_{i_1},\ldots,x_{i_s})y_{i_1\cdots i_t}(x_{i_1},\ldots,x_{i_t})f_{\bX}(\bx)\rm{d}\bx=0,
\end{equation}
where $i_1,\ldots,i_s \neq i_1,\ldots, i_t$ for $1\leq s \leq N$ and  $1\leq t \leq N$.
Applying the ANOVA dimensional decomposition in \eqref{ANOVA2} and the orthogonal properties in \eqref{ortho} to the variance of the quantity of interest, $\sigma_y=\int_{\mathbb{R}^N}(y(\bx)-y_0)^2 f_{\bX}(\bx)\rm{d}\bx$, results in the partitioning of the variance of $y(\bX)$. This is used for deriving the variance-based sensitivity method, as introduced in the following section.

\subsubsection{Variance-based sensitivity analysis}\label{sec:2.4:2}
For $\mathcal{V}\subseteq \{1,\ldots,N\}$, let the $|\mathcal{V}|$-variate global sensitivity index of a random output $Y = y(\bX)$ be 
\begin{equation}
S_{\mathcal{V}} := \dfrac{\sigma_{\mathcal{V}}^2}{\sigma_y^2}, \quad \sigma_y>0,
\label{eq:SV}
\end{equation}
where $\sigma_{\mathcal{V}}^2$ is the variance of $y_{i_1\ldots i_{|\mathcal{V}|}}(X_{i_1},\ldots,X_{i_{|\mathcal{V}|}})$ defined above. This non-negative sensitivity index reflects the fraction of the variance of $y(\bX)$ contributed by the inputs  $X_{i_1},\ldots,X_{i_{|\mathcal{V}|}}$.
%

\subsection{Polynomial dimensional decomposition surrogate}\label{sec:2.5}
%
\textcolor{black}{For the ANOVA dimensional decomposition of any output function $y(\cdot)$ that is square-integrable on the probability space $\left(\Omega,\mathcal{F},\mathbb{P}\right)$, the Cameron-Martin theorem~\citep{cameron1947orthogonal} states that there exists a corresponding Fourier-polynomial expansion.} This expansion is called \textit{polynomial dimensional decomposition (PDD)}, and is given by
\begin{align*}
y(\bX)=y_{0}+\sum_{\emptyset \neq \mathcal{U}\subseteq \{1,\ldots,N\}}\sum_{\bj_{\mathcal{U}}\in\mathbb{N}^{|\mathcal{U}|}}c_{\mathcal{U},\bj_{\mathcal{U}}}\Psi_{\mathcal{U},\bj_{\mathcal{U}}}(\bX_{\mathcal{U}}), \qquad \quad c_{\mathcal{U},\bj_{\mathcal{U}}}:=\int_{\mathbb{A}^N}y(\bx)\Psi_{\mathcal{U},\bj_{\mathcal{U}}}(\bx_{\mathcal{U}})f_{\bX}(\bx)\mathrm{d}\bx,
\end{align*}
where $c_{\mathcal{U},\bj_{\mathcal{U}}}$ are the expansion coefficients and the multivariate orthonormal polynomial is defined as 
$
\Psi_{\mathcal{U},\bj_{\mathcal{U}}}(\bX_{\mathcal{U}})=\prod_{i\in \mathcal{U}}\Psi_{i,j_i}(X_i),
$
where $\Psi_{i,j_i}$ is a univariate orthonormal polynomial in $X_i$ of degree $j_i$ that is consistent with the probability measure $f_{X_i}(X_i)\mathrm{d}x_i$. 
%
%
The full PDD contains an infinite number of expansion coefficients. In practice, the PDD must be truncated to have a finite number of expansion coefficients. A straightforward approach is to retain the degrees of interaction among input variables less than or equal to $S$. For example, when $S=1$ and $S=2$, the PDD includes at most univariate and bivariate polynomials, respectively. We then preserve the degree or order $m$ of the polynomial expansion such that $S \leq m < \infty$. This truncation results in an $S$-variate, $m$th-order \textit{PDD approximation}, i.e., 
\begin{align*}
y_{S,m}(\bX):=y_{0}+\sum_{\substack{\emptyset \neq \mathcal{U}\subseteq \{1,\ldots,N\}\\ 1\leq |\mathcal{U}| \leq S }}\sum_{\substack{\bj_{\mathcal{U}}\in\mathbb{N}^{|\mathcal{U}|}\\ |\mathcal{U}|\leq |\bj_{\mathcal{U}}|\leq m}}c_{\mathcal{U},\bj_{\mathcal{U}}}\Psi_{\mathcal{U},\bj_{\mathcal{U}}}(\bX_{\mathcal{U}})\approx y(\bX).
\end{align*} 
We can arrange the elements of the basis in any order, such that 
\begin{align*}
\{\Psi_{\mathcal{U},\bj_{\mathcal{U}}}(\bX_{\mathcal{U}})~:~1\leq |\mathcal{U}| \leq S,~|\mathcal{U}|\leq |\bj_{\mathcal{U}}| \leq m\}=\{\Psi_2(\bX),\ldots,\Psi_{L}(\bX)\},~\Psi_1(\bX)=1,
\end{align*}
where $\Psi_i(\bX)$ represents the $i$th basis function in the truncated PDD approximation. 
With this, the PDD approximation can be rewritten as
\begin{align}\label{pdd:3}
y_{S,m}(\bX)=\sum_{i=1}^{L}c_i\Psi_i(\bX)
\end{align}
where $c_i\in\mathbb{R}$ is the corresponding expansion coefficient for $i=1,\ldots,L$ and where
\begin{equation}  \label{eq:L}
L=L(S,N,m)=1+\sum_{s=1}^S\binom{N}{s}\binom{m}{s}.
\end{equation} 

We can use the $S$-variate, $m$th-order PDD approximation to estimate the global sensitivity index $S_{\mathcal{V}}$  from \eqref{eq:SV} for $\emptyset \neq \mathcal{V} \subseteq \mathcal{U}$ as 
\begin{align*}
S_{\mathcal{V}}\approx 
\left. 
\displaystyle\sum_{\substack{\bj_{\mathcal{V}}\in \mathbb{N}^{|\mathcal{V}|}\\|\mathcal{V}| \leq |\bj_{\mathcal{V}}|\leq m}}c_{\mathcal{V},\bj_{\mathcal{V}}}^2
\middle /
\displaystyle\sum_{\substack{\emptyset \neq \mathcal{U}\subseteq \{1,\ldots,N\}\\1\leq |\mathcal{U}| \leq S}}\sum_{\substack{\bj_{\mathcal{U}}\in\mathbb{N}^{|\mathcal{U}|}\\|\mathcal{U}| \leq |\bj_{\mathcal{U}}|\leq m}}c_{\mathcal{U},\bj_{\mathcal{U}}}^2. \right. 
\end{align*} 

\textcolor{black}{We refer the interested readers to ~\citep{sudret2008global}, which first introduced that the sensitivity indices can be computed directly from the PCE coefficients for the case when the input variables are uniformly distributed. The PDD surrogate is applicable to a wider range of distributions, including truncated normal, lognormal, and exponential distributions, yet the same formulas for sensitivity analysis apply.}

\subsection{Diffeomorphic Modulation under Observable Response Preserving Homotopy (D-MORPH) regression}\label{sec:2.6}
The authors in~\citep{li2010,li2012} introduce the D-MORPH regression to solve an underdetermined linear system, i.e., a system with more unknown parameters than training samples.
Consider $\bx^{(l)}=(x_1^{(l)},\ldots,x_N^{(l)})^{\intercal}$ for $l=1,\ldots, M$, obtained by (quasi) Monte Carlo or Latin hypercube sampling with corresponding probability $f_{\bX}(\bx)$. Given $M < L$, the expansion coefficients  $\bc=(c_1,c_2,\ldots,c_L)^{\intercal}\in\mathbb{R}^L$ in \eqref{pdd:3} of the PDD approximation can be obtained by solving an underdetermined linear system 
\begin{align}
\underbrace{
\begin{bmatrix}
\Psi_1(\bx^{(1)}) & \cdots &  \Psi_{L}(\bx^{(1)}) \\
\vdots                   & \ddots &  \vdots                           \\
\Psi_1(\bx^{(M)}) & \cdots &  \Psi_{L}(\bx^{(M)})
\end{bmatrix}}_{=:\bA } \underbrace{\begin{bmatrix} c_1 \\ c_2 \\ \vdots \\ c_L  \end{bmatrix}}_{\bc}
= \underbrace{\begin{bmatrix} y(\bx^{(1)}) \\ y(\bx^{(2)}) \\\vdots \\ y(\bx^{(M)}) \end{bmatrix}}_{=:\bb}.
\label{linear}
\end{align}
Assembling the right-hand side requires costly simulations, hence it is desirable to minimize the number of $M$. 
Since there exists an infinite number of solutions for $\bc$ that satisfy \eqref{linear}, a manifold $\mathcal{M}\subseteq \mathbb{R}^{L}$ is constructed to include all potential solutions. The D-MORPH regression aims to obtain the optimal solution $\bc$ within $\mathcal{M}$, by minimizing certain undesirable properties of $\bc$. 

Consider $\bA \in \mathbb{R}^{M \times L}$ with rank $r < L \leq M$ so that by the singular value decomposition
\begin{equation}
\bA =
\bH
\begin{bmatrix}
\bR_r      &  \boldsymbol{0} \\
\boldsymbol{0}    &  \boldsymbol{0}
\end{bmatrix} \bK^\intercal,
\end{equation} 
where, $\bH$ and $\bK$ are $M \times M$ and $L \times L$ orthogonal matrices, respectively, while $\bR_r$ is a nonsingular $r \times r$ diagonal matrix.  
The generalized inverse of the matrix $\bA$ is obtained as  
$
\bA^+=
\mathbf{K}
\small{
\begin{bmatrix}
\bR_r^{-1} &  \boldsymbol{0} \\
\boldsymbol{0}    &  \boldsymbol{0}
\end{bmatrix}} \bH^\intercal.
$
Here, $\bA^+ \in \mathbb{R}^{L \times L}$ satisfies the four Moore-Penrose conditions: $\bA \bA^+ \mathbf{A}=\bA$, $\bA^+ \bA \bA^+=\bA^+$, $(\bA \bA^+)^\intercal=\bA \bA^+$, $(\bA^+ \bA)^\intercal=\bA^+ \bA$. 

Consider $t\in \mathbb{R}_0^+$ as a scalar variable that is used to parametrize the vector $\ba(t)=(a_1(t),\ldots,a_{L}(t))^{\intercal}\in\mathbb{R}^{L}$ of all potential solutions of PDD expansion coefficients. 
Let $\mathbf{u}(t)$ denote an arbitrary function vector in $\mathbb{R}^{L}$. The set of potential solutions within $\cM$ is expressed as
\begin{align}
\ba(t) = \bA^+ \bb + (\bI_{L} - \bA^+ \bA)\bu(t).
\label{dmorph:1}
\end{align}
The first part of the solution $\bA^+ \bb$ in \eqref{dmorph:1} is an initial estimate for $\bc$, akin to the standard least-squares solution when the regularization term is ignored. 
We proceed to discuss the second term in~\eqref{dmorph:1}.
Define
\begin{align}
\boldsymbol{\Phi}:=(\bI_{L} - \bA^+ \bA) \in \mathbb{R}^{L \times L}\quad
\text{and}\quad
\bv(t):=
\displaystyle
\frac{\mathrm{d}\bu(t)}{\mathrm{d}t} \in \mathbb{R}^{L}.
\label{OG_proj}
\end{align}
From the Moore-Penrose conditions, we can show that $\boldsymbol{\Phi}$ is an orthogonal projector with the properties  $\boldsymbol{\Phi}^2=\boldsymbol{\Phi}$ and $\boldsymbol{\Phi}^\intercal = \boldsymbol{\Phi}$.
When differentiating \eqref{dmorph:1} with respect to $t$ and using \eqref{OG_proj}, we obtain 
\begin{align}\label{dmorph:ode}
\displaystyle
\frac{\mathrm{d}\ba(t)}{\mathrm{d}t}  = \boldsymbol{\Phi} \bv(t).
\end{align}
In the standard D-MORPH regression, one defines a quadratic cost function $\mathcal{K}(\ba(t)) \in \mathbb{R}_0^+$ and subsequently minimizes it, i.e., the original D-MORPH solves  
\begin{align}\label{costf2}
\min_{t\in\mathbb{R}}\left\{\mathcal{K}(\ba(t)) =
\displaystyle
\frac{1}{2}\ba^\intercal(t) \bD \ba(t)\right\}.
\end{align}
Here, $\mathbf{D}$ is an $L \times L$ real-valued, symmetric, non-negative definite matrix, so that the coefficients $a_i(t)$, $i=1,\ldots,L$, contract during the D-MORPH iterations at rates depending on the elements of $\mathbf{D}$.  If $\mathbf{D}$ is a diagonal matrix, then relatively larger values may be assigned to appropriate diagonal entries to suppress contributions from high-order basis functions of PDD.
In~\eqref{dmorph:ode}, we select 
\begin{align}\label{free_func}
\mathbf{v}(t) =
\displaystyle
-\frac{\partial \mathcal{K}(\ba(t))}{\partial \ba(t)}.
\end{align}
Using the chain rule and properties of the projector $\boldsymbol{\Phi}$, it can be shown that
\begin{align}
\displaystyle
\frac{\mathrm{d} \mathcal{K}(\ba(t))}{\mathrm{d}t} = -
\left( \boldsymbol{\Phi} \displaystyle \frac{\partial \mathcal{K}(\ba(t))}{\partial \ba(t)} \right)^\intercal
\left( \boldsymbol{\Phi} \displaystyle \frac{\partial \mathcal{K}(\ba(t))}{\partial \ba(t)} \right) \le 0.
\label{costf1}
\end{align}
According to \eqref{costf1}, the quadratic cost function $\mathcal{K}(\mathbf{a}(t))$ monotonically decreases as $t$ increases. 
Combining \eqref{dmorph:ode}, \eqref{costf2}, and  \eqref{free_func} results in an initial-value problem governed by the differential equation
\begin{align}\label{dmorph:2}
\displaystyle
\frac{\mathrm{d}\ba(t)}{\mathrm{d}t} = - \boldsymbol{\Phi} \bD \ba(t),\qquad
\ba(0) = \bA^+ \bb.
\end{align}
From \eqref{dmorph:2}, a transient solution is derived analytically, such that
\begin{equation}
\mathbf{a}(t) =
\exp(- t \boldsymbol{\Phi} \bD)  \ba(0) =
\exp(- t \boldsymbol{\Phi} \mathbf{D}) \bA^+ \bb.
\end{equation}
The singular value of decomposition of $\bPhi\bD$ is
\begin{align}\label{svd}
\boldsymbol{\Phi} \bD =
\bE
\begin{bmatrix}
\bT_r      &  \boldsymbol{0} \\
\boldsymbol{0}    &  \boldsymbol{0}
\end{bmatrix} \bF^\intercal,
\end{align}
with $\bT_r$ representing an $r \times r$ diagonal matrix with nonzero entries. By taking the limit $t \to \infty$, the final D-MORPH solution is 
\begin{align} 
\bar{\bc} =
\displaystyle \lim_{t \to \infty} \ba(t) =
\bF_{L-r} (\bE_{L-r}^\intercal \bF_{L-r})^{-1}
\bE_{L-r}^\intercal  ~\bA^+ \bb, 
\label{dmorph:sol}
\end{align}
where $\bar{\bc}=(\bar{c}_1,\ldots,\bar{c}_{L})^{\intercal}$ with $\bar{c}_i\in\mathbb{R},~i=1,\ldots,L$, representing the expansion coefficients of PDD by the original D-MORPH regression and matrices $\bE_{L-r}$ and $\bF_{L-r}$ are constructed from the last ${L-r}$ columns of matrices $\bE$ and $\bF$ from \eqref{svd}.


\section{A PDD surrogate obtained from a novel D-MORPH regression} \label{sec:3} 
We propose a new D-MORPH regression method to train the PDD surrogate in the limited data setting where we define a new cost function for the D-MORPH regression that includes a sparse Lasso solution. This method results in more robustness to variations in the inputs to the QoI while maintaining the training data fit. Section~\ref{sec:3.1} presents the challenges encountered when solving underdetermined systems for PDD. In Section~\ref{sec:3.2}, we propose the novel Lasso-based D-MORPH regression.  
\subsection{Challenges in solving underdetermined systems for PDD and Lasso regression}\label{sec:3.1} 
The standard PDD surrogate (detailed in Section~\ref{sec:2.5}) uses the ANOVA decomposition, and thus effectively truncates the basis functions that model higher-order interactions. The PDD surrogate thus mitigates the curse of dimensionality compared to the polynomial chaos expansion surrogate when the QoI is defined via a large number ($N$) of inputs or is heavily nonlinear so  that it requires a large degree ($m$) in the truncation of PDD. 
Nevertheless, the PDD surrogate can still require several hundreds to thousands of training samples. Acquiring such a large number of training samples can be computationally prohibitive. For example, assume that we require $500$ training samples for a bivariate ($S=2$) tenth-order ($m=10$) PDD approximation of five ($N=5$) random inputs to calculate the PDD expansion coefficients via a conventional regression method (e.g., standard least squares). If each sample requires a simulation that takes 24 CPU hours, then the process would take $12,000$ hours=$500\times24$ hours. For applications such as char combustion---the focus herein---we need to reduce the number of required training samples. 

We thus consider an underdetermined linear system for the expansion coefficients $\bc$ of an $S$-variate $m$th-order PDD approximation in~\eqref{linear}. 
When using the D-MORPH regression, the solution can be highly sensitive to the selection of the weight values in $\bD$ of \eqref{dmorph:2}. The strategies for selecting weights in \citep{li2010,li2012} are specific to a particular problem. In this work, we aim to generalize the strategy for selecting weights to work for different problems.   
Another option is to solve a Lasso regression 
\begin{align}
\min_{\bc\in\mathbb{R}^L}\left\{ (\bb-\bA\bc)^{\intercal}(\bb - \bA\bc) + k \sum_{i=1}^L|c_i|\right \}, 
\label{lasso}
\end{align}
where $k$ is a positive real number. The second term of \eqref{lasso} is
a regularization term that penalizes the $l^1$ norm of the PDD's expansion coefficients, producing sparse solutions for underdetermined systems. However, Lasso introduces a bias into the estimates to reduce the variance (the classical bias/variance trade-off in statistics). Due to the nature of the $l^1$ penalty, the Lasso cannot select more coefficients than the number of training samples. 

In the following section, we introduce a new D-MORPH regression method that combines the benefits of both the D-MORPH and the Lasso regression. The new method provides accurate PDD expansion coefficients from an underdetermined linear system~\eqref{linear}. This is possible as the new regression method improves the robustness of the D-MORPH solution to variations in input values or weight values and overcomes the disadvantages of the Lasso regression.

\subsection{Novel Lasso-based D-MORPH regression} \label{sec:3.2} 
To compute the PDD expansion coefficients $\bc$ for \eqref{pdd:3}, \textcolor{black}{ one typically solves the linear system~\eqref{linear}, which requires input-output data.} We consider the limited data setting, where the system for the unknown regression coefficients is underdetermined, i.e., $M < L$. \textcolor{black}{In the next sections, we define a new cost function and present an approach to solve the resulting optimization problem.}

\subsubsection{Cost function}\label{sec:3.2:1}
\textcolor{black}{We define a new cost function to minimize the $l^2$ norm of the difference between potential D-MORPH and Lasso solutions. Unlike the original D-MORPH cost function that only minimizes the solution's variance, the proposed cost function additionally promotes sparsity in the D-MORPH solution, similar to those obtained with Lasso. However, unlike Lasso which directly enforces sparsity, our approach allows coefficients to be near-zero values.} 

\textcolor{black}{Consider the case $M<L$, i.e., an underdetermined linear system~\eqref{linear} for the expansion coefficients of the PDD. Denote a D-MORPH solution to \eqref{linear} as $\ba(t)=(a_1(t),\ldots,a_{L}(t))^{\intercal}$, and denote the Lasso solution by $\bc_0=(c_{0,1},\ldots,c_{0,L})^{\intercal}$. In the initial iteration, we compute the starting D-MORPH solution as 
\begin{align}\label{prior_costf}
\bc_1=(c_{1,1},\ldots,c_{1,L})^{\intercal}=\argmin_{t\in\mathbb{R}}\left\{\dfrac{1}{2}(\mathbf{a}(t)-\mathbf{c}_0)^{\intercal}(\mathbf{a}(t)-\mathbf{c}_0)\right\}.
\end{align}}

\textcolor{black}{In the subsequent iteration, the cost function is augmented to include the $l^2$ norm between a potential D-MORPH solution and the solution obtained from the previous iteration.}
The new D-MORPH regression aims to minimize the cost function $\breve{\mathcal{K}}(\ba(t))$, i.e., 
\begin{align}\label{new_costf}
\min_{t\in\mathbb{R}}\left\{\breve{\mathcal{K}}(\mathbf{a}(t))=\dfrac{\lambda}{2}(\mathbf{a}(t)-\mathbf{c}_0)^{\intercal}\mathbf{W}(\mathbf{a}(t)-\mathbf{c}_0)+\dfrac{1-\lambda}{2}(\mathbf{a}(t)-\mathbf{c}_1)^{\intercal}\mathbf{W}(\mathbf{a}(t)-\mathbf{c}_1)\right\}
\end{align}
with a non-negative real-valued weight $\lambda\in[0,1]$. Here, $\bW=\rm{diag}[1/(c_{1,1}+\epsilon),\ldots,1/(c_{1,L}+\epsilon)]$ is an $L$-dimensional diagonal matrix, and where $\epsilon \ll 1$ (e.g., $\epsilon = 1\times 10^{-6}$). \textcolor{black}{This choice of weights aims to preserve the sparsity in $\mathbf{c}_1$. For example,} when the elements of $\bc_1$ are smaller, this matrix assigns larger weights to the corresponding elements of $\ba(t)-\bc_0$ and $\ba(t)-\bc_1$. 
As the zero entries of $\bc_1$ have the largest weights, $1/\epsilon$, the corresponding elements of $\ba(t)$ are the most strongly constrained to the elements of $\bc_0$ and $\bc_1$.    
\textcolor{black}{We note that the new cost function~\eqref{new_costf} decreases monotonically as $t$ increases during the D-MORPH process, as demonstrated in \eqref{costf1}.}   

\subsubsection{Non-homogeneous ordinary differential equation for D-MORPH regression}\label{sec:3.2:2}
We perform a D-MORPH regression to minimize the new cost function~\eqref{new_costf}. This way, we combine the advantages of both the D-MORPH and Lasso regression and obtain a more robust solution to the expansion coefficients of the PDD surrogate.
With the new cost function~\eqref{new_costf}, we set up a non-homogeneous ordinary differential equation in a similar fashion to \eqref{dmorph:ode} and \eqref{free_func} as  
\begin{align}
\dfrac{\rm{d}\ba(t)}{\rm{d}t}=-\bPhi \bW \ba(t)+\bPhi\bW \left (\lambda \bc_0 + (1-\lambda)\bc_1\right ),\qquad \ba(0)=\bA^{+}\bb, 
\label{nonhode}
\end{align}
which has the solution 
\begin{align}\label{dmorph:ode2}
\ba(t)=\exp(-t\mathbf{\Phi}\bW)\bA^+\bb+\displaystyle\int_0^t\exp\{-(t-q)\mathbf{\Phi}\bW\}\mathbf{\Phi}\bW({\lambda\bc_0}+(1-\lambda)\bc_1)\mathrm{d}q,
\end{align}
where the second term represents the particular solution.
Taking the limit $t\rightarrow\infty$, the D-MORPH solution to \eqref{linear} is $\breve{\bc}=\lim_{t\rightarrow \infty}\ba(t)$, which can be written as
\begin{equation}
\begin{split} 
\breve{\bc}=\bar{\bF}_{L-r}(\bar{\bE}^{\intercal}_{L-r}\bar{\bF}_{L-r})^{-1}\bar{\bE}^{\intercal}_{L-r}\bA^+\bb
+\bar{\bF}_{r}(\bar{\bE}^{\intercal}_{r}\bar{\bF}_{r})^{-1}\bar{\bE}^{\intercal}_{r}(\bar{\bT}_r)^{-1}\mathbf{\Phi}\bW^{-1}(\lambda\bc_0+(1-\lambda){\bc_1}).
\end{split}
\label{dmorph:sol2}
\end{equation}
Here, $\bar{\mathbf{E}}_r$ and $\bar{\mathbf{E}}_{L-r}$, $\bar{\mathbf{F}}_r$ and $\bar{\mathbf{F}}_{L-r}$ are constructed from the first $r$ and the last $L-r$ columns of matrices $\bar{\mathbf{E}}$ and $\bar{\mathbf{F}}$, respectively, generated from the singular value decomposition
\begin{equation}
\mathbf{\Phi}\bW=\bar{\bE}
\begin{bmatrix}
\bar{\bT}_r & \mathbf{0}\\
\mathbf{0} & \mathbf{0}
\end{bmatrix}
\bar{\bF}^{\intercal},
\end{equation}
with $\bPhi$ representing an orthogonal projector in \eqref{OG_proj} and $\bar{\bT}_r$ representing an $r\times r$ diagonal matrix including nonzero singular values. 
%
\subsubsection{Recursive process for improving the D-MORPH solution}\label{sec:3.2:3}
\textcolor{black}{We now formulate the new iterative version of the D-MORPH solution $\breve{\bc}$. In \eqref{dmorph:sol2}, we replace the best D-MORPH solution $\breve{\bc}$ \eqref{dmorph:sol2} with the D-MORPH solution  $\breve{\bc}^{(i)}=(\breve{c}_1^{(i)},\ldots,{\breve{c}_{L}}^{(i)})^{\intercal}$ at iteration $i$. To ensure convergence of the iteration, we redefine the prior solution $\bc_1$ as  
\begin{align} 
\bc_1^{(i-1)} = \dfrac{1}{i}\sum_{k=1}^{i} \breve{\bc}^{(k-1)}. 
\end{align} 
This formula represents the average of all prior solutions computed up to the $(i-1)$th iteration. These refinements yield}
\begin{equation} 
\begin{split}
\breve{\bc}^{(i)}=&\bar{\bF}^{(i-1)}_{L-r}(\bar{\bE}^{(i-1)\intercal}_{L-r}\bar{\bF}^{(i-1)}_{L-r})^{-1}\bar{\bE}^{(i-1)\intercal}_{L-r}\bA^+\bb+
\\&\bar{\bF}^{(i-1)}_{r}(\bar{\bE}^{(i-1)\intercal}_{r}\bar{\bF}^{(i-1)}_{r})^{-1}\bar{\bE}^{(i-1)\intercal}_{r}(\bar{\bT}_r^{(i-1)})^{-1}\mathbf{\Phi}\mathrm{diag}({\breve{\bc}^{(i-1)}})^{-1}(\lambda\bc_0+(1-\lambda){\bc_1}^{(i-1)}),
\end{split}
\end{equation}
where $\bar{\mathbf{E}}^{(i)}_r$ and $\bar{\mathbf{E}}^{(i)}_{L-r}$, $\bar{\mathbf{F}}^{(i)}_r$ and $\bar{\mathbf{F}}^{(i)}_{L-r}$ are constructed from the first $r$ and the last $L-r$ columns of matrices $\bar{\mathbf{E}}^{(i)}$ and $\bar{\mathbf{F}}^{(i)}$, respectively, generated from the singular value decomposition
\begin{equation}
\mathbf{\Phi}\bW^{(i-1)}=\bar{\mathbf{E}}^{(i-1)}
\begin{bmatrix}
\bar{\mathbf{T}}_r^{(i-1)} & \mathbf{0}\\
\mathbf{0} & \mathbf{0}
\end{bmatrix}
\bar{\mathbf{F}}^{(i-1)\intercal},
\end{equation}
with $\bar{\mathbf{T}}^{(i-1)}$ representing an $r\times r$ diagonal matrix including nonzero singular values. Here, $\bW^{(i-1)}=\rm{diag}[0,1/(\breve{c}_2^{(i-1)}+\epsilon),\ldots,1/(\breve{c}_{L}^{(i-1)}+\epsilon)]$, where its first element is \emph{zero} and the remaining elements are the reciprocal of  $\breve{c}_j^{(i-1)}$ for $j=2,\ldots,L$ with $\epsilon \ll 1$. 

The initial D-MORPH regression is 
\begin{equation}
\breve{\bc}^{(0)}=\bar{\mathbf{F}}_{L-r}(\bar{\mathbf{E}}^{\intercal}_{L-r}\bar{\mathbf{F}}_{L-r})^{-1}\bar{\mathbf{E}}^{\intercal}_{L-r}\mathbf{A}^+\mathbf{b}
+
\bar{\mathbf{F}}_{r}(\bar{\mathbf{E}}^{\intercal}_{r}\bar{\mathbf{F}}_{r})^{-1}\bar{\mathbf{E}}^{\intercal}_{r}(\bar{\mathbf{T}}_r)^{-1}\mathbf{\Phi}\breve{\bc}_0,
\end{equation} 
where $\bar{\mathbf{E}}_r$ and $\bar{\mathbf{E}}_{L-r}$, $\bar{\mathbf{F}}_r$ and $\bar{\mathbf{F}}_{L-r}$ are constructed from the first $r$ and the last $L-r$ columns of matrices $\bar{\mathbf{E}}$ and $\bar{\mathbf{F}}$, respectively. These matrices are generated from the singular value decomposition
\begin{equation}
\mathbf{\Phi}=\bar{\mathbf{E}}
\begin{bmatrix}
\bar{\mathbf{T}}_r & \mathbf{0}\\
\mathbf{0} & \mathbf{0}
\end{bmatrix}
\bar{\mathbf{F}}^{\intercal},
\end{equation}
with $\bar{\mathbf{T}}_r$ representing an $r\times r$ diagonal matrix including nonzero singular values.

\section{Global sensitivity analysis by a PDD surrogate}
\label{sec:4} 
We leverage the new D-MORPH regression PDD surrogate for global sensitivity analysis in the limited data setting.
Given a limited computational budget, Section~\ref{sec:4.1} presents the complete algorithm of the proposed method for global sensitivity analysis. 
In Section~\ref{sec:4.2}, we evaluate the proposed method for global sensitivity analysis using the Ishigami \& Homma function. We assess the convergence of the proposed method by comparing its results to the exact sensitivity solutions using the Oakley \& O'Hagan function in Section~\ref{sec:4.3}. 

\subsection{Complete algorithm for global sensitivity analysis}\label{sec:4.1}
\begin{figure*}[ht]
\begin{center}
\includegraphics[angle=0,scale=0.90,clip]{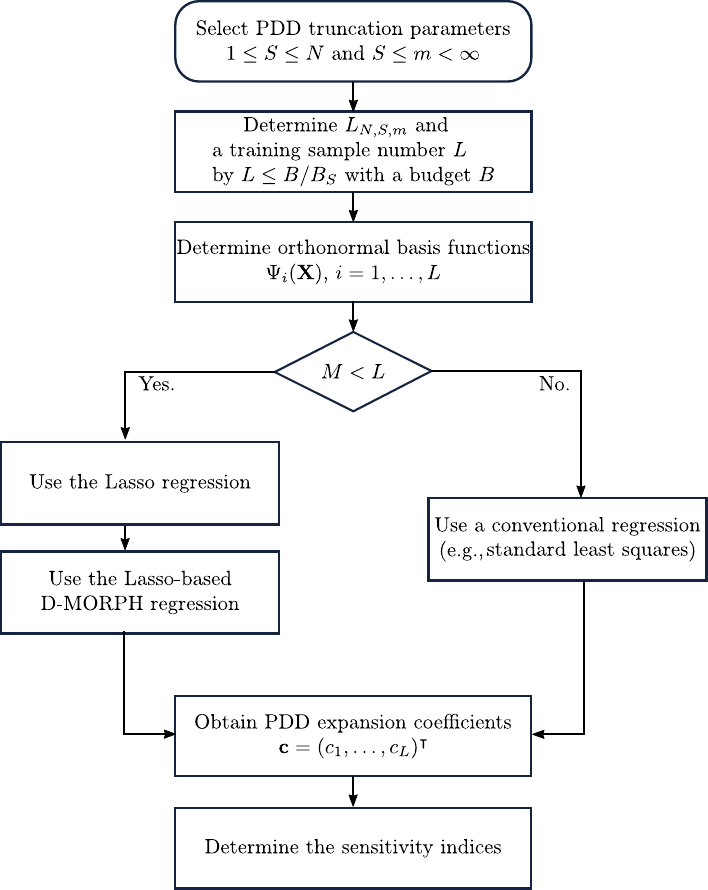}
\end{center}
\caption{Flow chart for computing the PDD surrogate with the proposed D-MORPH regression for global sensitivity analysis.}
\label{fig1}
\end{figure*} 
The flow chart in Figure~\ref{fig1} details the procedure for implementing the PDD surrogate modeling for global sensitivity analysis under a budget constraint. 
Define $B\in\mathbb{R}_0^+$ as the total computational budget and $B_S\in\mathbb{R}_0^+$ as the computational cost required to obtain a single training sample. The total cost is $B=M\times B_S$, where $M$ is the number of training samples affordable by the computational budget.  
We then select the PDD truncation parameters $S$ and $m$. In this work, we choose $S=2$ and $m$ ranges from 5 to 11 to attain a target estimate that aligns closely with an unbiased Monte Carlo estimate with $M$ samples.       
In the limited data setting, $M<L$,  we use Lasso regression to obtain the initial $\bc_0$ in \eqref{new_costf} followed by the proposed D-MORPH regression to obtain PDD expansion coefficients~\eqref{dmorph:sol2}. Otherwise, we have the overdetermined case and can use a conventional regression (e.g., standard least squares). 
We determine the first and the second sensitivity indices and the total effect sensitivity indices, as detailed in Section~\ref{sec:2.5}. 

\subsection{Illustrative example~1: Ishigami \& Homma function}\label{sec:4.2}
This first illustrative example demonstrates that the iterative process improves the accuracy of the D-MORPH solution for the highly nonlinear output with high-order interactions.
\subsubsection{Problem definition}\label{sec:4.2:1} 
Consider the Ishigami \& Homma function from \citep{ishigami1990importance}, given as 
\begin{align} 
y(\mathbf{X})=\sin{X_1}+a\sin^2{X_2}+bX_3^4\sin{X_1},
\end{align} 
where $X_1, X_2, X_3$ are independent and identically distributed uniform input random variables on $[-\pi,+\pi]$, and $a$ and $b$ are real-valued deterministic parameters; we select $a=7$ and $b=0.1$. Table~\ref{table1} reports the exact solutions for the variance and sensitivity indices of the output random variable $Y=y(\bX)$. 
\begin{table*}
\caption{Exact solutions for the mean, the standard deviation, and the first-order and the second-order sensitivity indices of Ishigami \& Homma function.}
\vspace{-0.1in}
\begin{spacing}{1.2}
\begin{center}
\small
\resizebox{\textwidth}{!}{
\begin{tabular}{ccc}
\toprule 
 & Exact solution & Exact solution ($a=7$, $b=0.1)$\tabularnewline
\midrule
\textcolor{black}{$\mu_{Y}$ (Mean of $Y$)} & \textcolor{black}{$a/2$} & \textcolor{black}{$3.5$}\tabularnewline
$\sigma_{Y}$ (Standard deviation of $Y$) & $\sqrt{a^{2}/8+b\pi^{4}/5+b^{2}\pi^{8}/18+1/2}$ & $3.720832$$^{\mathrm{(a)}}$\tabularnewline
$S_{\{1\}}$ (First-order sobol index for $X_{1}$) & $(b\pi^{4}/5+b^{2}\pi^{8}/50+1/2)/\sigma_{Y}^{2}$ & $0.313905$$^{\mathrm{(a)}}$\tabularnewline
$S_{\{2\}}$ (First-order sobol index for $X_{2}$) & $(a^{2}/8)/\sigma_{Y}^{2}$ & $0.442411$$^{\mathrm{(a)}}$\tabularnewline
$S_{\{3\}}$ (First-order sobol index for $X_{3}$) & $0$ & $0$\tabularnewline
$S_{\{1,2\}}$ (Second-order sobol index for $X_{1}$ and $X_{2}$) & $0$ & $0$\tabularnewline
$S_{\{1,3\}}$ (Second-order sobol index for $X_{1}$ and $X_{3}$) & $(\ensuremath{8b^{2}\pi^{8}}/\ensuremath{225})/\ensuremath{\sigma_{Y}^{2}}$ & $0.243684$$^{\mathrm{(a)}}$\tabularnewline
$S_{\{2,3\}}$ (Second-order sobol index for $X_{2}$ and $X_{2}$) & $0$ & $0$\tabularnewline
\bottomrule
\end{tabular}}
\end{center}
\end{spacing}
\vspace{-0.15in}
\begin{tablenotes}
\scriptsize\smallskip 
\item{a.} The exact solution is rounded to six decimal places. 
\end{tablenotes}
\label{table1}
\end{table*}
\subsubsection{Process of global sensitivity analysis}
We select the PDD truncation parameters $S=2$ and $m=11$ so that the PDD has $L=199$ expansion coefficients, see \eqref{eq:L}. 
We consider $M=59$ training samples, which is $30\%$ of $L=199$. Since the linear system from \eqref{linear} is underdetermined (i.e., $M < L$), we use the Lasso-based D-MORPH regression and select the weight $\lambda=0.5$. 
We then  use the obtained expansion coefficients to determine the variance and the first and second-order sensitivity indices, as detailed in Sections~\ref{sec:2.4} and \ref{sec:2.5}. 
\subsubsection{Results}
\begin{table*}
\caption{Standard deviation and first-order and second-order sensitivity estimates by bivariate ($S=2$) $m=11$th-order PDD approximation using Lasso-based DMORPH regression and Lasso regression, with 59 training samples (30\% of the number of unknown expansion coefficients).}
\vspace{-0.1in}
\begin{center}
\small
\resizebox{\textwidth}{!}{
\begin{tabular}{clccccccc}
\toprule 
\multicolumn{2}{c}{} & \multicolumn{4}{c}{Mean relative error$^\mathrm{a}$ ($K=30$)} & \multicolumn{3}{c}{Mean absolute error$^{\mathrm{b}}$ ($K=30$)}\tabularnewline
\multicolumn{2}{c}{Methods} & Standard deviation & $S_{\{1\}}$ & $S_{\{2\}}$ & $S_{\{1,3\}}$ & $S_{\{3\}}$ & $S_{\{1,2\}}$ & $S_{\{2,3\}}$\tabularnewline
\midrule
\multicolumn{2}{l}{Lasso-based D-MORPH} &  &  &  &  &  &  & \tabularnewline
 & Iteration=$0$ & $0.128965$ & $0.057175$ & $0.103271$ & $0.192922$ & $0.002249$ & $0.002876$ & $0.008772$\tabularnewline
 & Iteration=$20$ & $0.083872$ & $0.035874$ & $0.073809$ & $0.149792$ & $0.001271$ & $0.001382$ & $0.005195$\tabularnewline
 & Iteration=$30$ & $0.069427$ & $0.028573$ & $0.078167$ & $0.137125$ & $0.000573$ & $0.001464$ & $0.004164$\tabularnewline
\multicolumn{2}{l}{Lasso regression} & $0.157845$ & $0.064819$ & $0.112890$ & $0.199169$ & $0.001420$ & $0.000820$ & $0.005834$\tabularnewline
\bottomrule
\end{tabular}}
\end{center}
\vspace{-0.15in}
\begin{tablenotes}
\scriptsize\smallskip 
\item{a.} The mean relative error over 30 trials is the mean absolute error over 30 trials, normalized by the exact solution.
\item{b.} The mean absolute error over 30 trials is used when the exact solution is \emph{zero}.  
\end{tablenotes}
\label{table2}
\end{table*}
Table~\ref{table2} reports the mean errors (absolute and relative, depending on if we have zero or non-zero reference values) of the bivariate eleventh-order PDD computed by the Lasso-based D-MORPH regression for the standard deviation, the first-order sensitivity indices $S_{\{i\}}$, and the second-order sensitivity indices $S_{\{i,j\}}$, where $i,j=1,2,3$ and $j>i$, compared to the exact solutions in Table~\ref{table1}. The \textit{mean absolute error} and the \textit{mean relative error} are obtained by
\begin{align*}
\dfrac{1}{K}\sum_{k=1}^K\left|\mathcal{Y}-\mathcal{Y}_{k}\right|~\qquad\text{and}\qquad \dfrac{1}{K}\sum_{k=1}^K\left|\dfrac{\mathcal{Y}-\mathcal{Y}_{k}}{\mathcal{Y}}\right|,
\end{align*}
where $\mathcal{Y}$ is an exact solution for the standard deviation or the sensitivities of $y(\bX)$ and $\mathcal{Y}_k$ is the PDD estimate of $\mathcal{Y}$ at the $k$th independent trial run. When the exact solution is non-zero, we normalize the mean absolute error with respect to the corresponding exact solution and obtain the mean relative error.
As shown in the second through the fourth row of Table~\ref{table2}, increasing the iteration number $i$ of Equation~\eqref{dmorph:sol2} from 0 to 30 decreases the mean errors accordingly. In particular, the mean relative error of the proposed D-MOPRH for the standard deviation decreases by almost 50\% from iteration 0 to 30. 

We also include the Lasso estimates in the last row of Table~\ref{table2}. It shows that in the standard deviation and most cases of sensitivity indices, the proposed D-MORPH estimates are more accurate than the Lasso estimates, while they use only 30\% of the number of expansion coefficients. 
\begin{figure}[htbp]
\centering 
\begin{subfigure}{0.496\textwidth}{\includegraphics[width=\linewidth, trim={0.4cm 0.4cm 2cm 3cm}, clip]{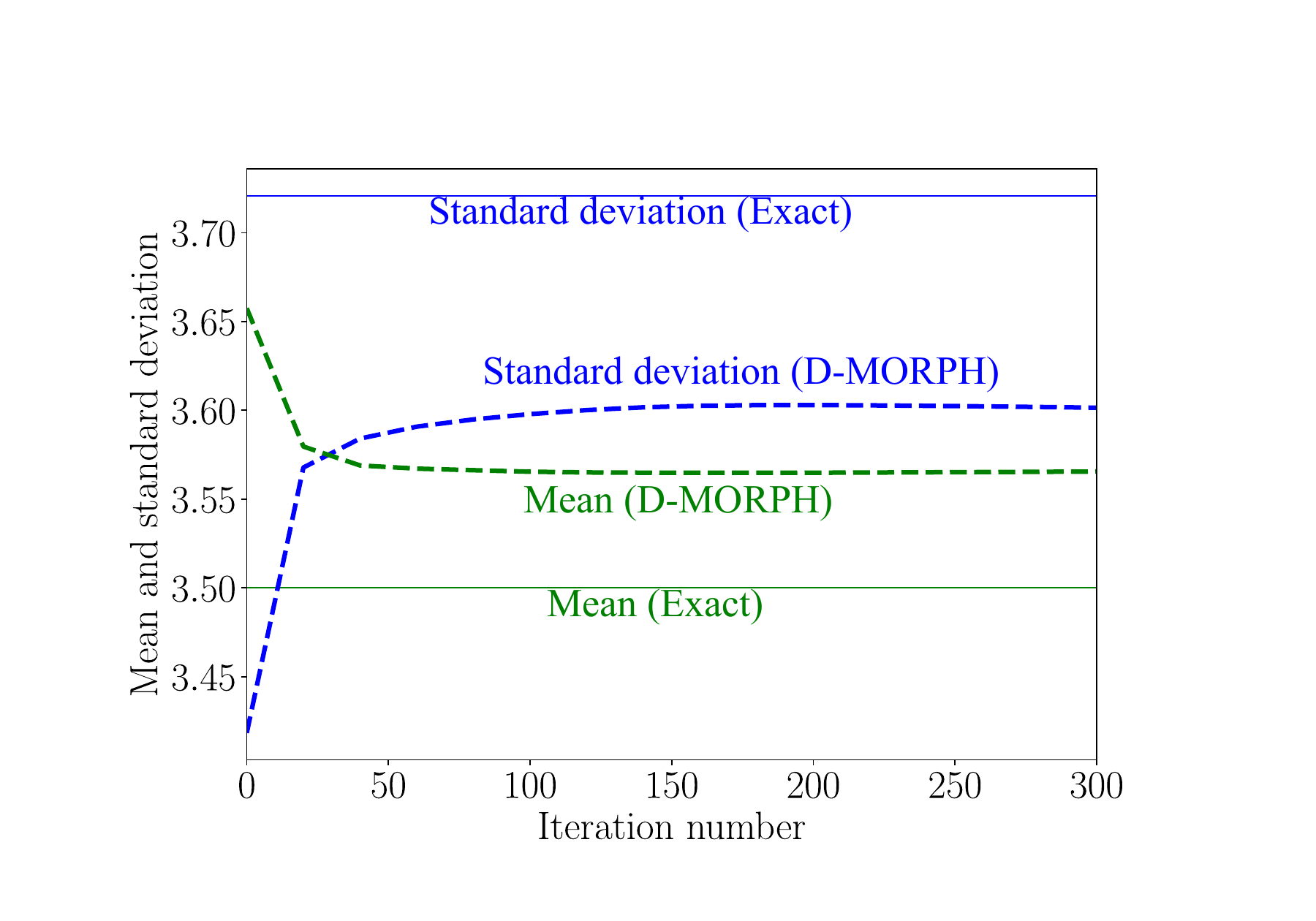}}\caption{\textcolor{black}{Mean and standard deviation}}\label{fig2:(a)}
\end{subfigure} 
\begin{subfigure}{0.496\textwidth}{\includegraphics[width=\linewidth, trim={0.4cm 0.4cm 2cm 3cm}, clip]{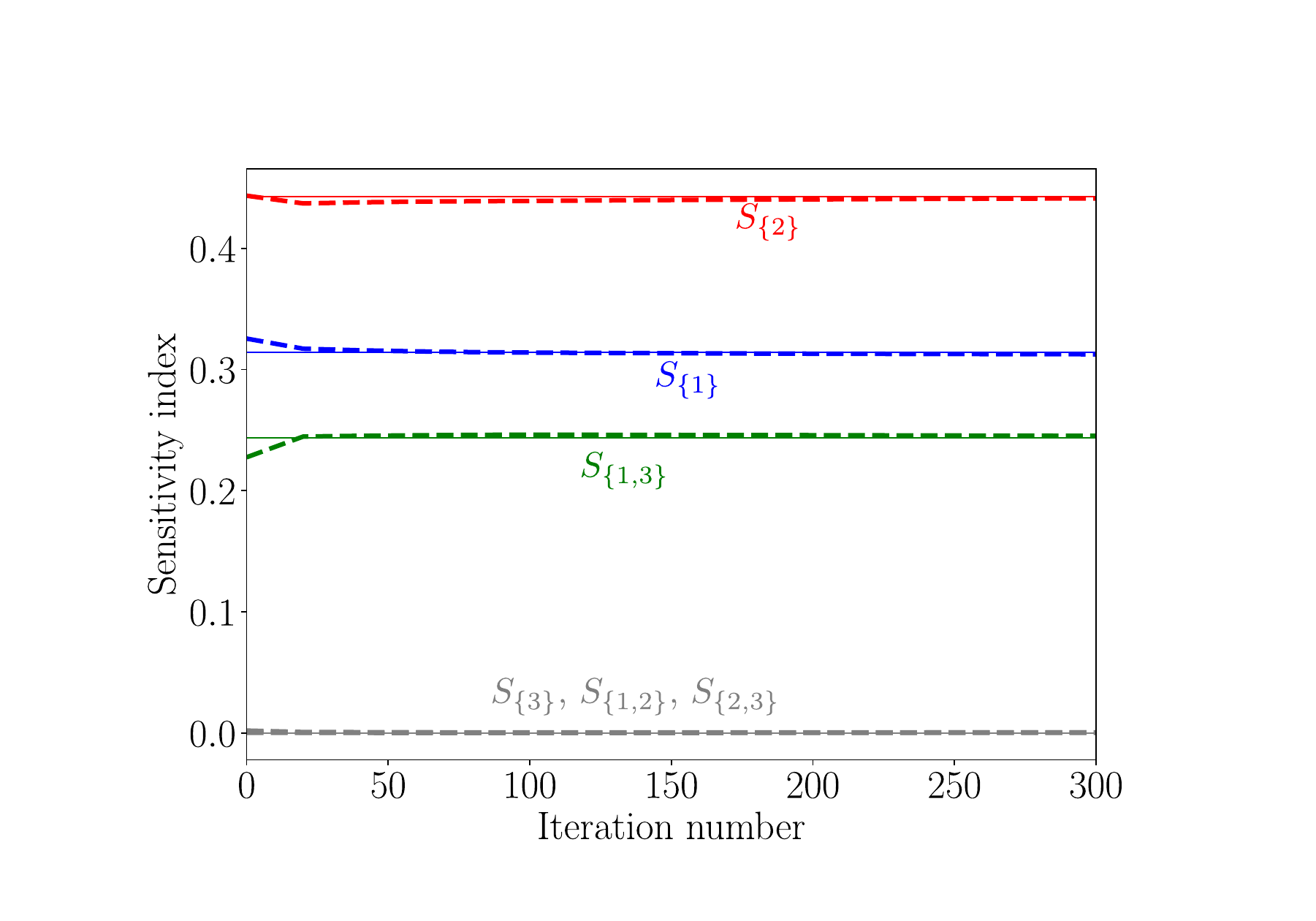}}\caption{Sensitivity indices}\label{fig2:(b)}
\end{subfigure} 
\begin{subfigure}{0.496\textwidth}
\end{subfigure} 
\label{fig:2}
\caption{(a) Convergence of PDD estimates computed by a Lasso-based D-MORPH regression \textcolor{black}{for mean} and standard deviation; (b) sensitivity indices $S_{\{1\}}, S_{\{2\}}, S_{\{3\}}$ as D-MORPH iteration increases from 0 to 300.}
\end{figure} 
Figures~\ref{fig2:(a)} and \ref{fig2:(b)} show the convergence of one realization of PDD estimates computed by the proposed D-MOPRH regression as the iteration number ($i$) of \eqref{dmorph:sol2} increases from 0 to 300. In Figure~\ref {fig2:(a)}, the mean and the standard deviation converge relatively rapidly, yet there exists a bias. For example, within 40 iterations, \textcolor{black}{the mean and standard deviation estimates converge from 4.57\% to 2.9\% error and from 8.18\%  to 3.66\% error, respectively, when compared to their respective exact solutions.
On the other hand,} the sensitivities converge rapidly to each of the exact solutions yet without significant bias. For example, within 40 iterations, the first-order sensitivity estimate $S_{\{1\}}$ converges from 3.69\%  to 0.61\% error when compared to the exact solution. After 40 iterations, the convergence rate for the standard deviation and sensitivities ($S_{\{1\}}$, $S_{\{3\}}$, $S_{\{1,3\}}$) tends to slow down.

\subsection{Illustrative example 2: Oakley \& O'Hagan function}\label{sec:4.3}  
In this second example for a relatively high-dimensional global sensitivity analysis problem, we demonstrate the convergence of D-MORPH regression as the number of training samples increases. We also consider different weight values ($\lambda$) for the D-MORPH cost function~\eqref{new_costf}.

\subsubsection{Problem definition}
For global sensitivity analysis, Oakley \& O'Hagan~\citep{oakley2004probabilistic} introduced a mixture of trigonometric and quadratic polynomial functions as 
\begin{align*}
y(\mathbf{X})=\mathbf{a}_1^{\intercal}\mathbf{X}+\mathbf{a}_2^{\intercal}\sin\mathbf{X}+\mathbf{a}_3^{\intercal}\cos\mathbf{X}+\mathbf{X}^{\intercal}\mathbf{M}\mathbf{X}^{\intercal},
\end{align*}
where $\bX=(X_1,\ldots, X_{15})^{\intercal}\in \mathbb{R}^{15}$ is the standard Gaussian input vector ($N=15$) with mean vector $\mathbb{E}[\bX]=(0,\ldots,0)^{\intercal}\in\mathbb{R}^{15}$ and covariance matrix $\mathbb{E}[\bX\mathbf{X}^{\intercal}]-\bI\in\mathbb{R}^{15 \times 15}$. Moreover, $\mathbf{a}_i\in\mathbb{R}^{15}$, $i=1,2,3$, and $\bM\in\mathbb{R}^{15\times 15}$ are coefficient vectors and matrix, respectively, obtained from \citep{oakley2004probabilistic}. From the same work, we also obtain the exact solutions for the first-order sensitivities of the fifteen inputs.

\subsubsection{Process of global sensitivity analysis}
We choose the PDD truncation parameters $S=2$ and $m=5$. This results in a PDD with $L=1,126$ expansion coefficients. 
We consider three distinct training sample numbers $M=337$, $563$, and $788$, which are 30\%, 50\%, and 70\% of the unknown expansion coefficients. Since the linear system in \eqref{linear} is underdetermined (i.e., $M< L$), we use the Lasso-based D-MORPH regression. For the regression, we select three  weights $\lambda$=$0.2$, $0.6$, and $1.0$. 
\subsubsection{Results}\label{sec:3.5:2}
\begin{figure}[htbp]
\centering 
\includegraphics[width=0.8\linewidth, trim={8cm 2cm 8cm 6cm}, clip]{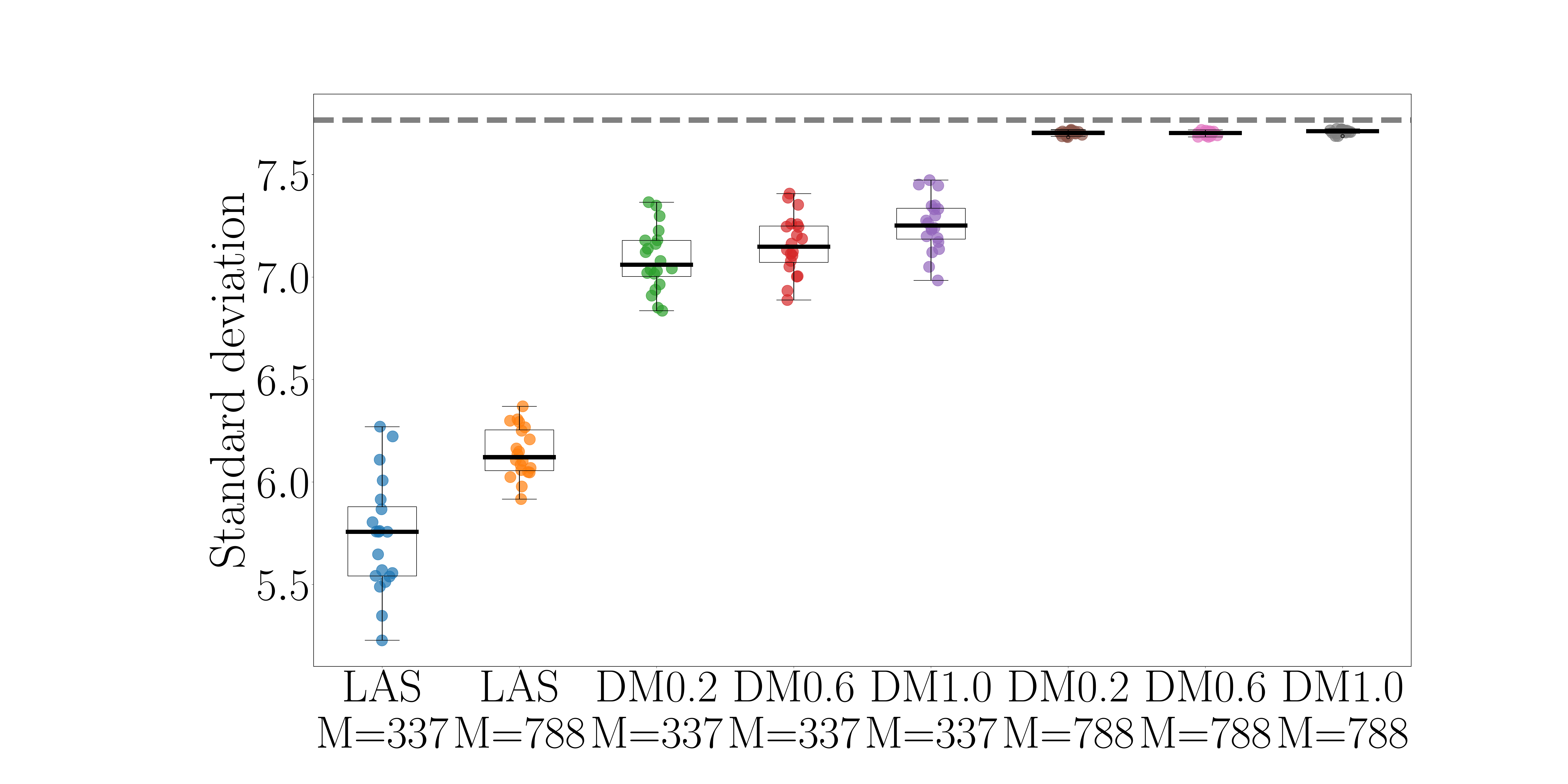}
\caption{Boxplots of the standard deviation of the random output $y(\bX)$, estimated by the bivariate fifth-order PDD using the Lasso (LAS) and Lasso-based D-MOPRH regressions with weight values $\lambda=0.2$, $0.6$, $1.0$ (DM0.2, DM0.6, DM1.0). Two underdetermined systems are considered: $M=337$ and $M=788$, which correspond to 30\% and 70\% of the number ($L =1,126$) of expansion coefficients. Each regression is repeated 20 times.  The exact solution is shown as a gray-dotted line.}
\label{fig3}
\end{figure} 
Figure~\ref{fig3} shows the standard deviations estimated by the bivariate fifth-order PDD approximations via box plots. These three distinct cases are associated with underdetermined linear systems from~\eqref{linear}. Hence, we use a Lasso regression and the proposed Lasso-based D-MORPH regression with three distinct weight values $\lambda=0.2$, $0.6$, and $1.0$, as indicated as `LAS, `DM0.2', `DM0.6', and `DM1.0', respectively, on the $x$-axis of Figure~\ref{fig3}. In the figure, we present the exact solution for the standard deviation with a gray dotted line.

\begin{figure}[h]
\centering 
\begin{subfigure}{0.495\textwidth}{\includegraphics[width=\linewidth, trim={8cm 2cm 8cm 6cm}, clip]{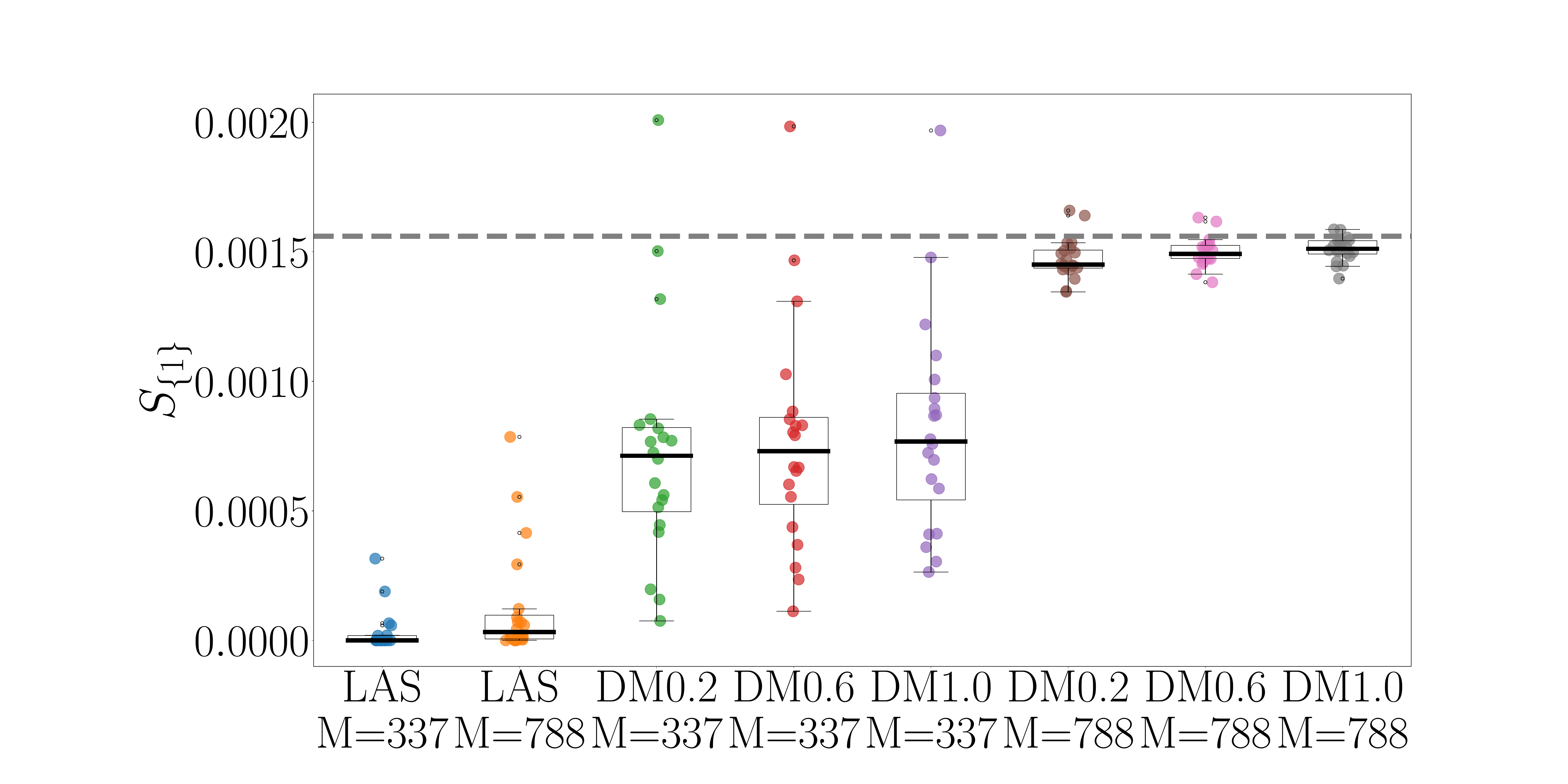}}\caption{$S_{\{1\}}$}\label{fig4:(a)}
\end{subfigure} 
\begin{subfigure}{0.495\textwidth}{\includegraphics[width=\linewidth, trim={8cm 2cm 8cm 6cm}, clip]{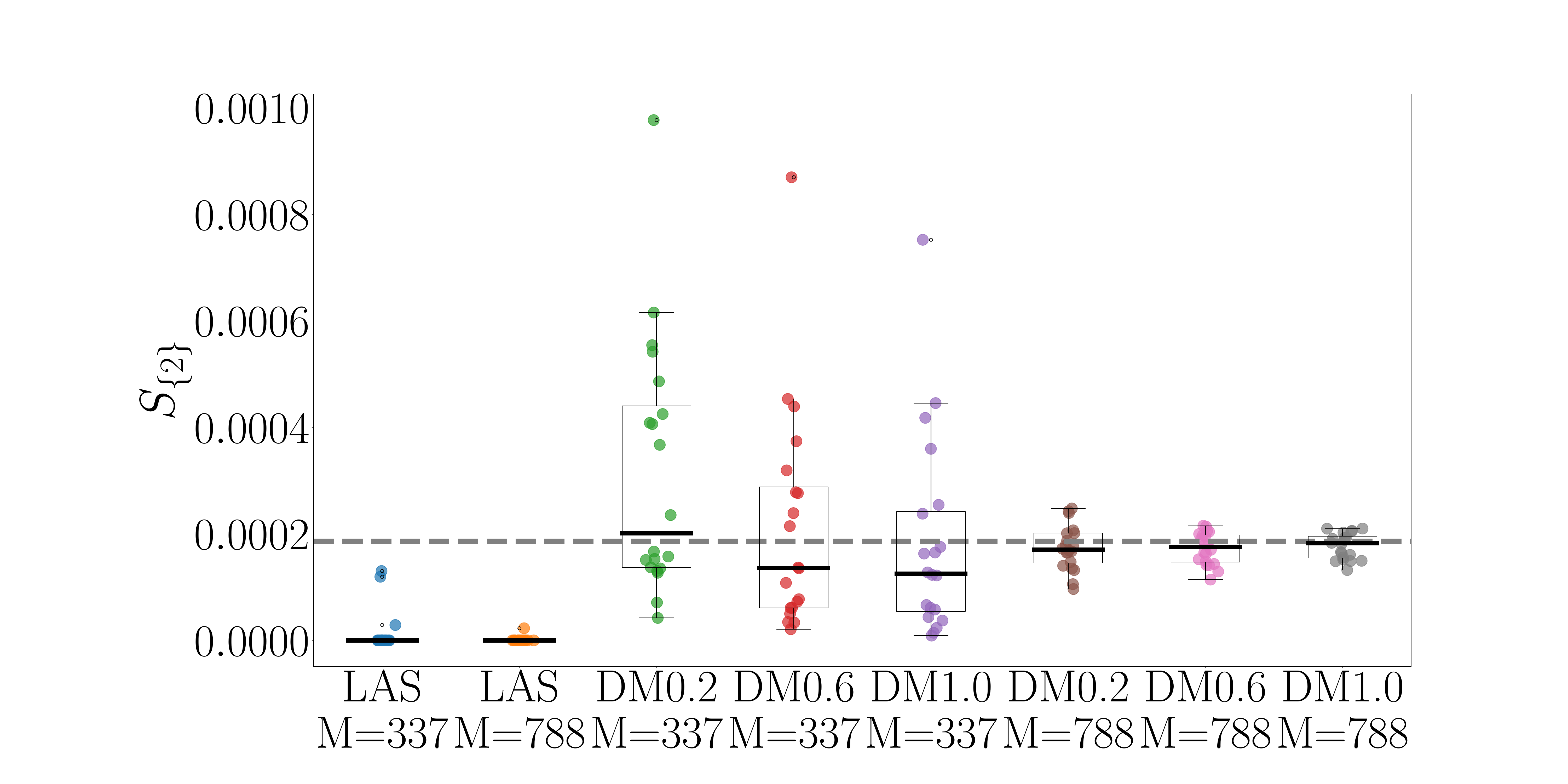}}\caption{$S_{\{2\}}$}\label{fig4:(b)}
\end{subfigure} 
\begin{subfigure}{0.495\textwidth}{\includegraphics[width=\linewidth, trim={8cm 2cm 8cm 5cm}, clip]{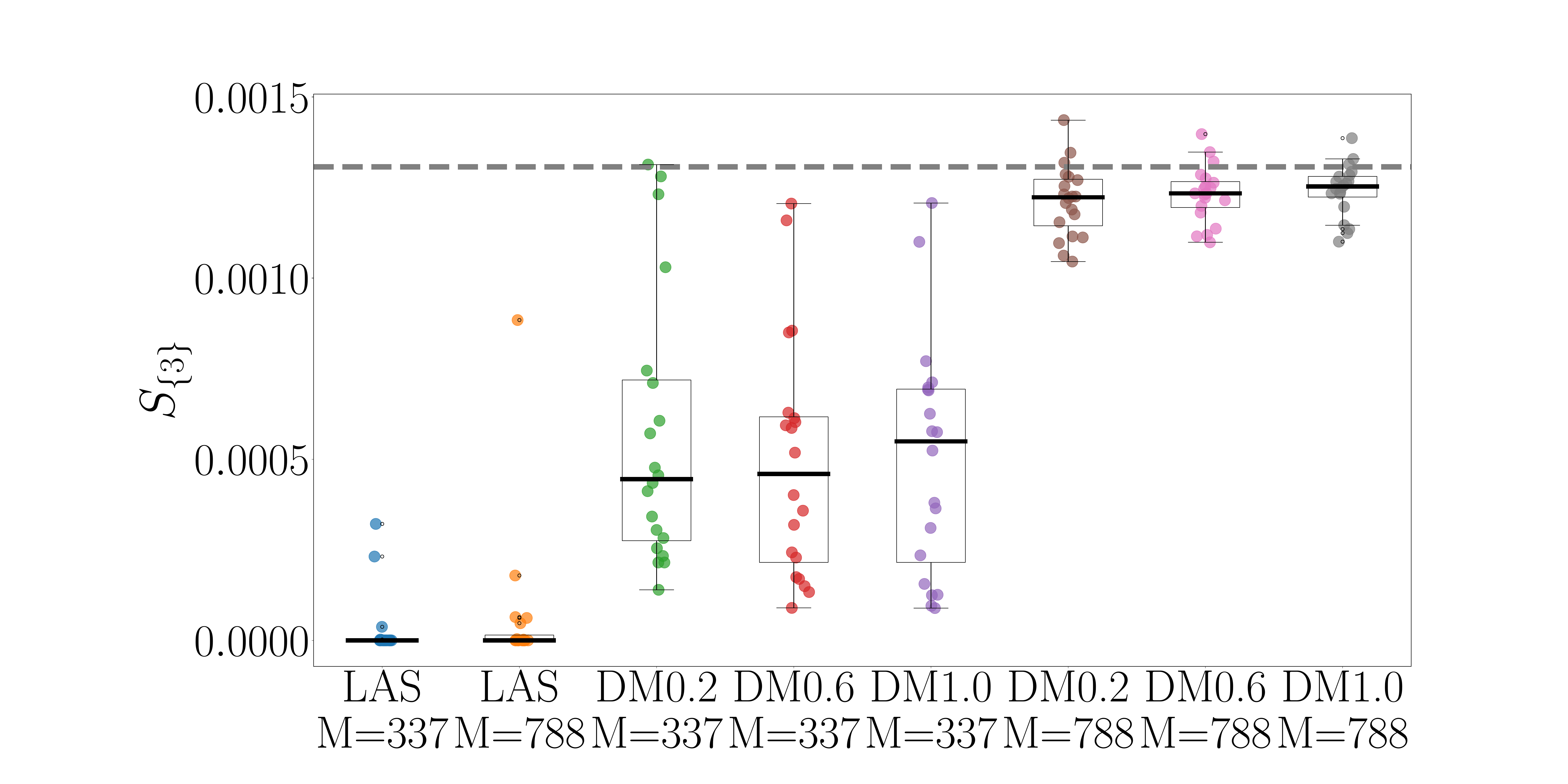}}\caption{$S_{\{3\}}$}\label{fig4:(c)}
\end{subfigure} 
\begin{subfigure}{0.495\textwidth}{\includegraphics[width=\linewidth, trim={8cm 2cm 8cm 5cm}, clip]{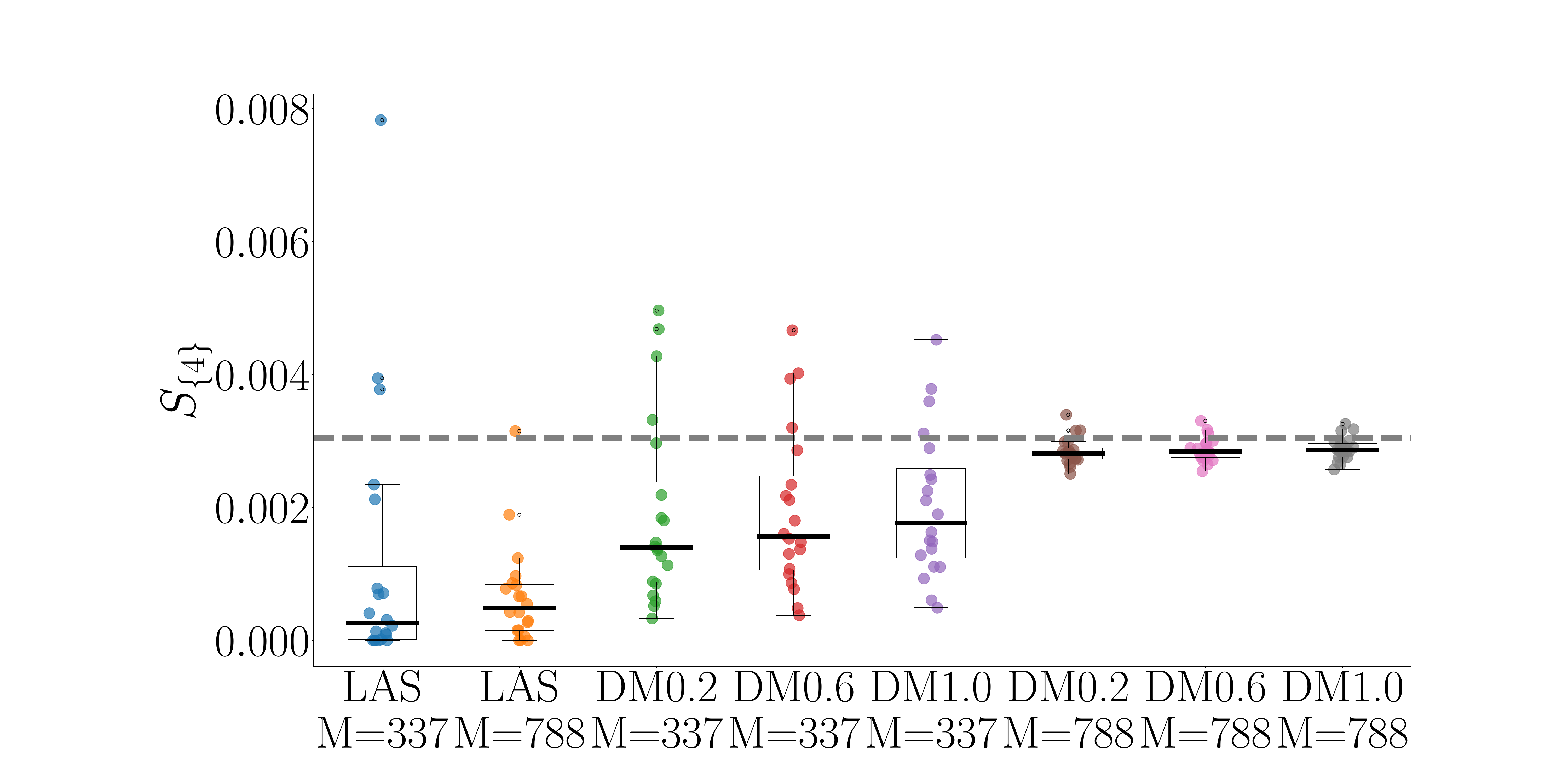}}\caption{$S_{\{4\}}$}\label{fig4:(d)}
\end{subfigure} 
\begin{subfigure}{0.495\textwidth}{\includegraphics[width=\linewidth, trim={8cm 2cm 8cm 5.5cm}, clip]{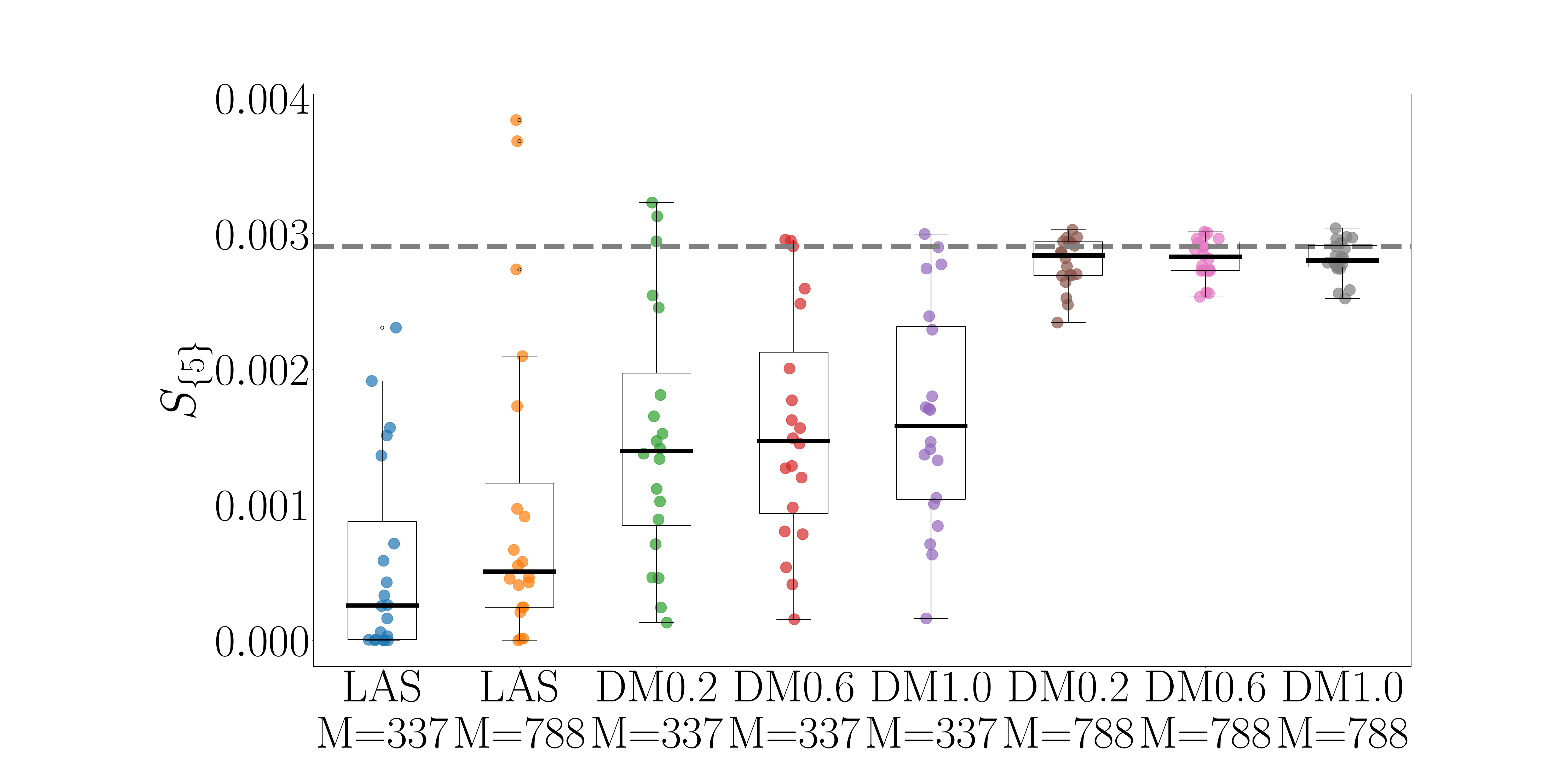}}\caption{$S_{\{5\}}$}\label{fig4:(e)}
\end{subfigure}
\begin{subfigure}{0.495\textwidth}{\includegraphics[width=\linewidth, trim={8cm 2cm 8cm 5.5cm}, clip]{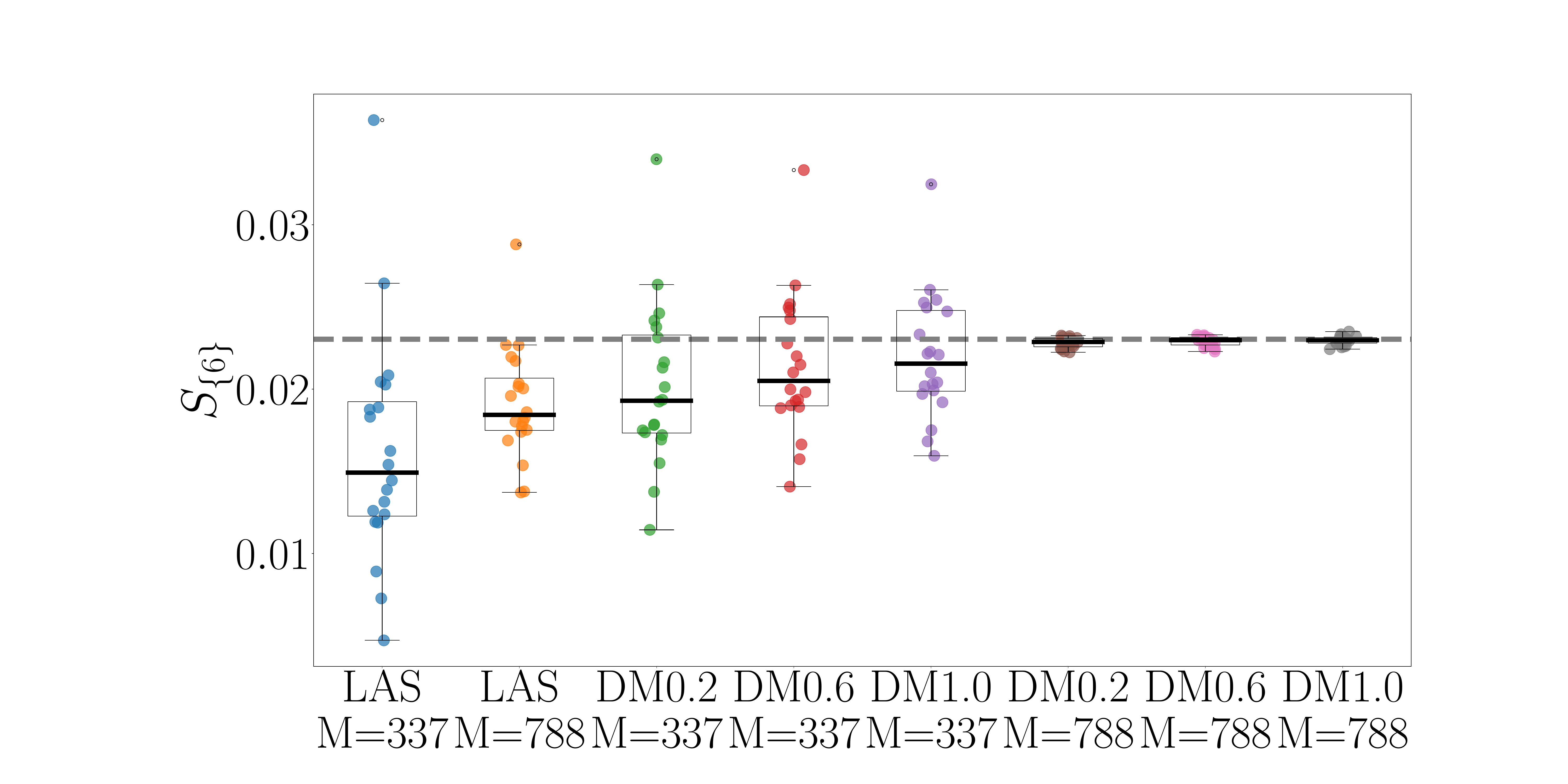}}\caption{$S_{\{6\}}$}\label{fig4:(f)}
\end{subfigure}
\caption{Boxplots of first-order sensitivity indices $S_{\{i\}}$, $i=1,2,3,4,5,6$, of the random output $y(\bX)$, in (a), (b), (c), (d), (e), (f), respectively, estimated by the bivariate fifth-order PDD using the Lasso (LAS) and Lasso-based D-MOPRH regressions with weight values $\lambda=0.2$, $0.6$, $1.0$ (DM0.2, DM0.6, DM1.0). Two underdetermined systems are considered: $M=337$ and $M=788$, which correspond to 30\% and 70\% of the number ($L=1,126$) of expansion coefficients. Each regression is repeated 20 times.  The exact solution is shown as a gray-dotted line.}
\label{fig4}
\end{figure}

As the number of training samples increases from 337 to 788, the proposed D-MORPH regressions for all weight cases converge more closely to the exact solution compared to the Lasso regression. For example, when $M=337$, the mean relative errors of the D-MORPH-based estimates compared to the exact standard deviation are 8.74\%, 7.84\%, and 6.56\% for weights 0.2, 0.6, and 1.0, respectively, over 20 experiments. These are 3x more accurate than the Lasso estimates (the mean relative error 26.18\% over 20 experiments). 
The D-MORPH regressions are more robust than the Lasso regression, as shown by the smaller variance of estimates in 20 experiments. As the number of training samples increases, the variance of the D-MORPH regressions also becomes narrower. \textcolor{black}{Figures~\ref{fig4:(a)}--\ref{fig4:(d)} show the exact values of $S_{\{i\}}$, for $i=1,2,3,4$, as dash-dotted lines. These values are nearly zero. In these sensitivity cases, LAS predicts a value of zero, whereas D-MORPH predicts a non-zero value with a wider variance over 20 experiments at a smaller sample size ($M=337$). As the sample number $M$ increases to 788 (70\% of the unknown expansion coefficients), the D-MORPH estimates are more accurate and show a narrower variance over 20 experiments compared to the LAS.}    

Figure~\ref{fig4} shows the estimates obtained by Lasso and the proposed D-MORPH regression for the first-order sensitivity indices $S_{\{i\}}$, $i=1-6$, which are representative of the other indices as well.  In Appendix~\ref{sec:appx1}, Figures~\ref{fig10} and \ref{fig11} show the results for the first-order sensitivity indices $S_{\{i\}}$, $i=7-15$. Taken together, these results demonstrate the superiority of the proposed D-MORPH method for three weight values ($\lambda=0.2$, $0.6$, $1.0$) over the Lasso regression in the first-order sensitivity indices, akin to the standard deviation case discussed earlier. Relative to number of training samples $L$, the estimates of the sensitivity indices computed from the D-MORPH regression are less sensitive to the weight $\lambda$. We note that $\lambda=1$ shows slightly higher accuracy compared to the other weight values.  

\section{Numerical example: Char combustion} \label{sec:5} 
Combustion is a computationally expensive process to simulate. Training surrogate models with such expensive simulation data presents challenges, as one is limited by how much data can be generated with a realistic computational budget. 
The proposed regression method addresses this computational challenge. In this section, we evaluate the PDD computed by the proposed D-MORPH regression for global sensitivity analysis of an expensive-to-simulate char combustion process with five random input variables.  

Sections~\ref{sec:5.1} and \ref{sec:5.2} describe the details of the problem and its numerical setting. We validate the simulation model in Section~\ref{sec:5.3}. Section~\ref{sec:5.4} clarifies the quantity of interest for global sensitivity analysis. We present the proposed surrogate modeling results for the mean and the standard deviation of the quantity of interest in Section~\ref{sec:5.5}. Finally, Section~\ref{sec:5.6} shows the results for sensitivity indices via the proposed methods. 
%

\subsection{Problem description}\label{sec:5.1}
\begin{figure*}
\begin{center}
\begin{subfigure}{0.43\textwidth}
{\includegraphics[width=\linewidth]
{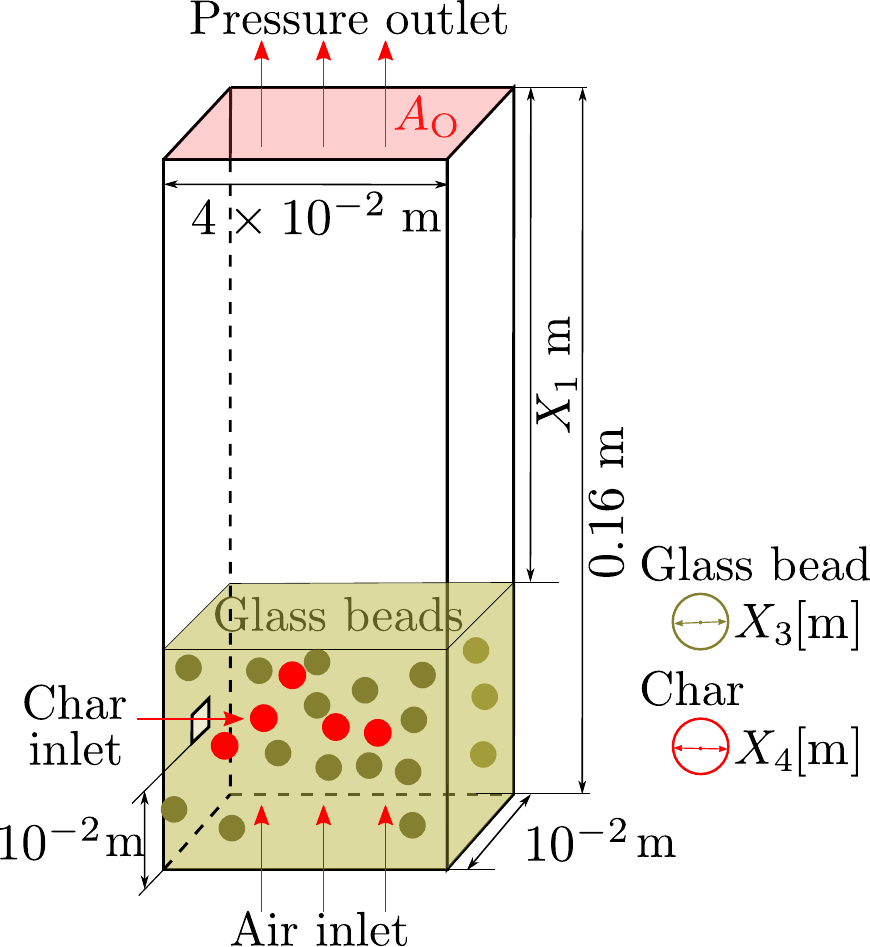}}\caption{Fluidized bed for char combustion}\label{fig5a}
\end{subfigure} 
\begin{subfigure}{0.27\textwidth}
{\includegraphics[width=\linewidth]
{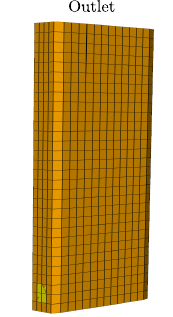}}\caption{Cell model for gas}\label{fig5b}
\end{subfigure} 
\hspace{0.15in}
\begin{subfigure}{0.25\textwidth}
{\includegraphics[width=\linewidth]
{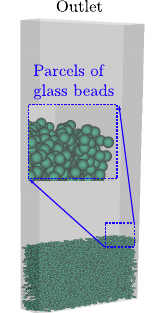}}\caption{Parcel model for glass beads}\label{fig5c}
\end{subfigure} 
\end{center}
\caption{Fluidized bed for char combustion: (a) the schematic diagram shows the geometry and the initial concentration of glass beads with a diameter of $X_3$ stacked at the boiler's bottom to a height of $X_1$. Char with a diameter of $X_4$ is fed in through the left side; (b) the cell model used to predict gas behavior consists of 2,520 cells; and (c) the parcel model used to predict solid behavior contains 32,945 parcels.}
\label{fig5}
\end{figure*} 

Fluidized bed combustion is a combustion technology that burns solid fuels, such as char and biomass, efficiently and with low emissions. The fluidized bed combustion systems can also capture pollutants, such as sulfur dioxide and nitrogen oxides, making them a more environmentally friendly option than traditional combustion methods.
To optimize the operational efficiency of fluidized bed combustors, this study focuses on determining the influential parameters that affect the QoI, here a thermal energy, which is computed from combustion simulations. We perform variance-based global sensitivity analysis using the PDD surrogate model with the proposed D-MORPH regression.  

Figure~\ref{fig5a} shows a geometrical configuration of a fluidized bed for char combustion~\citep{xie2021study}. The rectangular boiler is a lab-scale model with dimensions of $0.04$ m $\times$ $0.16$ m $\times$ $0.01$ m in width, height, and thickness. 
In the boiler model, we consider a total of five ($N=5$) random inputs, so $\bX = (X_1, X_2, X_3, X_4, X_5)^{\intercal}$. \textcolor{black}{The boiler initially contains a freeboard of height $X_1$ [m] that consists only of gases. Below this, there is a layer of glass beads with a particle diameter of $X_4$ [m].} The char particles react with oxygen from the air at a constant rate of $X_2$, generating heat and other products as a char combustion process. Char particles with a diameter of $X_3$ [m] are fed into the boiler at a constant rate of $X_5$ through the $2$ m $\times$ $2$ m char inlet on the left wall.   Table~\ref{table3} lists the random inputs, their interval bounds, and assumptions on their distribution.
\begin{table*}
\caption{Properties of the random inputs in a fluidized bed model for char combustion.}
\vspace{-0.1in}
\begin{center}
\small
\resizebox{\textwidth}{!}{
\begin{tabular}{ccccccc}
\toprule 
Random  & \multirow{2}{*}{Property} & \multirow{2}{*}{Mean} & \multirow{2}{*}{COV (\%)} & Lower  & Upper & Probability \tabularnewline
variable &  &  &  & boundary & boundary & distribution\tabularnewline
\midrule
$X_{1}$ & \textcolor{black}{Height of freeboard (m)} & $0.4$ & - & $\textcolor{black}{0.10}$ & $\textcolor{black}{0.14}$ & Uniform\tabularnewline
$X_{2}$ & Air inflow (m/s) & $0.825$ & $10$ & $0.425$ & $1.225$ & Truncated normal\tabularnewline
$X_{3}$ & \textcolor{black}{Diameter of the char particle (m)} & $8\times10^{-4}$ & - & $2\times10^{-4}$ & $1.4\times10^{-3}$ & Uniform\tabularnewline
$X_{4}$ & \textcolor{black}{Diameter of the glass bead particle (m)} & $1\times10^{-3}$ & - & $\color{black}{5\times10^{-5}}$ & $1.5\times10^{-3}$ & Uniform\tabularnewline
$X_{5}$ & Char mass inflow (kg/s) & $7.35\times10^{-6}$ & $10$ & $1.35\times10^{-6}$ & $1.35\times10^{-5}$ & Truncated normal\tabularnewline\midrule
\end{tabular}}
\end{center}
\vspace{-0.15in}
\label{table3}
\end{table*}
\subsection{Numerical setting} \label{sec:5.2}
To predict the combustion behavior in the fluidized bed boiler, we create a numerical model using particle-in-cell (PIC), which we summarize in this section.
\subsubsection{Particle-in-cell} \label{sec:5.2:1}
In the PIC, the gas phase is modeled using the Eulerian method, which treats the gas phase as continua. Figure~\ref{fig5b} shows the computational boiler  model used for gas phase simulation, consisting of 2,520 cells. Each cell has six degrees of freedom associated with three velocity components and three scalar variables (temperature, species concentrations, and pressure).    

The solids phases of the PIC are modeled using the Lagrangian method, where  particles with the same physical properties (e.g., density and diameter) are grouped to effectively track their positions and trajectories. A group of particles with the same physical properties is called a parcel.  Different particle diameters with the same materials are viewed as having different physical properties; therefore, they will have separate parcels. Figure~\ref{fig5c} shows the parcel model composed of 8,344 parcels for representing the glass beads in the boiler model.    
Since the PIC method tracks parcels and not particles, PIC provides a significant computational cost reduction. Nevertheless, the combustion simulation is expensive since it incorporates complex particle-scale physics such as chemical reactions, heat transfer, hydrodynamics, and more.
The PIC uses the collisional stress model to account for interactions between particles and walls~\citep{snider2001incompressible}. The collision stress is expressed as
\begin{align*}
\tau = \frac{P_{s}\varepsilon_{s}^{\beta}}
{\max(\varepsilon_{cp}-\varepsilon_{s},\alpha(1-\varepsilon_{s}))},
\end{align*} 
where $\varepsilon_{cp}$ is the pre-determined value that describes the maximum possible packing fraction for particles. Here, $P_{s}$ and $\beta$ are an empirical pressure constant and an empirical unit-less exponent, and $\alpha=10^{-9}$ is a non-singularity constant $\alpha\ll 1$. In this work, we select $P_{s}=1.0$ and $\beta=2.0$. 
We use the open-source software MFiX (version 23.1.1)~\citep{clarke2020mfix} with an MPI-based parallel computing solver on 15 CPUs (Intel Xeon W-3175X CPU @ 3.10 GHz) for the combustion simulations.

\subsubsection{Governing equations for particle-in-cell} \label{sec:5.1:1}
In this section, we briefly summarize the governing equations applied in the PIC method to predict the behaviors of both the gas and solid phases. We then introduce the chemical reaction used to predict the char combustion process.

\paragraph{Gas phase model} 
The conservation of mass, momentum, and internal energy are formulated~\citep{musser2020theoretical} as 
\begin{align*}
\frac{\partial}{\partial t}(\varepsilon_{g}\rho_{g}) +  \frac{\partial}{\partial x_j} (\varepsilon_{g}\rho_{g}U_{gj}) &= \sum_{n=1}^{N_g}R_{gn}  + S_g,\\
\frac{\partial }{\partial t}(\varepsilon_{g}\rho_{g}U_{gi}) +  \frac{\partial}{\partial x_j} (\varepsilon_{g}\rho_{g}U_{gj}U_{gi}) &= -\frac{\partial P_g}{\partial x_i} + \frac{\partial \tau_{gij}}{\partial x_j} + \varepsilon_{g}\rho_{g}g_{i} + S_{gi}, \\
\varepsilon_{g}\rho_{g}C_{pg}\left[\frac{\partial T_g}{\partial t} +  U_{gj}\frac{\partial T_g}{\partial x_j}\right] &= -\frac{\partial}{\partial x_j}(\varepsilon_{g}q_{gj}) + S_{g}.
\end{align*}
Here, in the gas phase, $\varepsilon_{\text{g}}$ is the volume fraction [$\rm{dimensionless}$], $\rho_{g}$ is the density [$\rm{kg/m^3}$], $P_g$ is the pressure~[$\rm{Pa}$], $T_g$ is the temperature~[$\rm{K}$], and $C_{pg}$ is the mixture of specific heat [$\rm{J/(kg\cdot K)}$]; $S_{g}$ is a general user-defined source~[$\rm{kg/(m^3\cdot s)}$] and $S_{gi}$ is a momentum term~[$\rm{kg/(m^2\cdot s^2)}$]. Additionally, $N_g$ is the number of chemical species in the gas phase and $R_{gn}$ is the rate of formation~[$\rm{kmole/ (s\cdot m^{3})}$] of the $n$th gas phase, while $g_i$ is the gravitational force~[$\rm{N}$], $U_{ji}$ is the velocity of gas~[$\rm{m/s}$] in the $j$th direction, and $\tau_{gij}$ is the stress tensor~[$\rm{Pa}$] in the gas phase. 

\paragraph{Solid phase model} 
The conservation of mass, momentum, and internal energy are formulated \citep{clarke2020mfix} as  
\begin{align*}
\frac{\partial}{\partial t}(W_{p}m_p)  &= W_p \sum_{n=1}^{N_p} R_{pn},\\
\dfrac{\partial}{\partial t}(W_{p}m_{p}X_{pn}) &= W_pR_{pn}, \\
W_{p}m_{p}\frac{\partial U_i}{\partial t} &= W_p(m_{p}g_{i} + \frac{m_p}{\varepsilon_{s}\rho_{s}} \nabla_{\vec{x}} \tau_p), \\
W_{p}m_{p}C_{p}\frac{\partial T}{\partial t} &= -W_p \sum_{n=1}^{N_p} h_{pn}R_{pn} + S_{p}.
\end{align*}
Here, $W_p$ is the statistical weight~[$\rm{dimensionless}$] of the particle; for the $n$th chemical species, $R_{pn}$ is the rate of the production/consumption~[$\rm{kmole/ s}$], $X_{pn}$ is the mass fraction~[$\rm{dimensionaless}$], and $h_{pn}$ is the specific enthalpy~[$\rm{J/kg}$]. Additionally, $C_p$ is the specific heat~[$\rm{J}/(kg\cdot K)$] and $T$ is the temperature~[$\rm{K}$] of the parcel, while $S_{p}$ is the general source term~[$\rm{kg/(m^3\cdot s)}$] of the parcel.

\paragraph{Chemical reaction}
Char combustion is governed by a heterogeneous gas-solid chemistry. The combustion process generates the two gas products $\rm{CO}$ and $\rm{CO}_2$ from
\begin{equation*}
\rm{C}(\text{Solid}) + 0.5\rm{O_2}(\text{Gas}) \rightarrow \rm{CO}(\text{Gas})
\quad\text{and}\quad 
\rm{CO}(\text{Gas}) + 0.5\rm{O_2}(\text{Gas}) \rightarrow \rm{CO_2}(\text{Gas}).
\end{equation*}
The collision among particles and the wall leads to the ash falling off from the particles, following an Arhenius kinetic rate and gas diffusion rate \citep{xie2021study}, i.e.,

\begin{equation*}
\frac{{\rm{d}}m_{ci}}{\rm{d}t} = -\pi d^{2}p_{o}\bigg(\frac{1}{R_{\rm{diff}}} + \frac{1}{R_{\rm{chem}}}\bigg)^{-1},~ R_{\rm{diff}} = \frac{24ShD_{o}}{d_{p}RT_{m}},~ R_{\rm{chem}} = A_i\exp\bigg(-\frac{E_i}{RT_p}\bigg),~ d_p = \bigg(\frac{6m_p}{\pi\rho_{p}}\bigg)^{1/3}.
\end{equation*}

Here, $m_{ci}$ is the unreacted char mass~[$\rm{kg}$]; $R_{\rm{diff}}$ is the gas diffusion rate~[$\rm{m^2/s}$] and $R_{\rm{chem}}$ is the Arrhenius kinetic rate~[$\rm{s^{-1}}$]. Also, $Sh$ is the Sherwood number~[$\rm{dimensionless}$], $R$ is the gas constant~[$\rm{J/(mol\cdot K)}$], and $T$ is the temperature~[$\rm{K}$]. Additionally, $p_o$ is the oxygen partial pressure~[$\rm{Pa}$] and $D_o$ is the oxygen-nitrogen mixture diffusion coefficient~[$\rm{m^2/s}$], while $d_p$ is the particle size~[$\rm{m}$] as it shrinks due to the mass loss. 

The homogeneous reaction can be calculated by the law of mass action via the Arrhenius formula, proposed by Dryer and Glassman \citep{dryer1973high}, i.e., 
\begin{equation*}
r_{co} = 3.98\times 10^{14} \exp\bigg(-\frac{1.67\times 10^5}{RT_g}\bigg)C_{\rm{CO}}C_{\rm{O_2}}^{0.25}C_{\rm{H_2O}}^{0.5},
\end{equation*}
where $C_{\rm{CO}}$, $C_{\rm{O}_2}$, and $C_{\rm{H_2O}}$ are the mass concentrations [$\rm{kmole/m^3}$] of $\rm{CO}$, $\rm{O}_2$, and $\rm{H_2O}$, respectively.   

\subsection{Validation} \label{sec:5.3}
\begin{figure*}
\begin{center}
\includegraphics[angle=0,scale=1.0,clip]{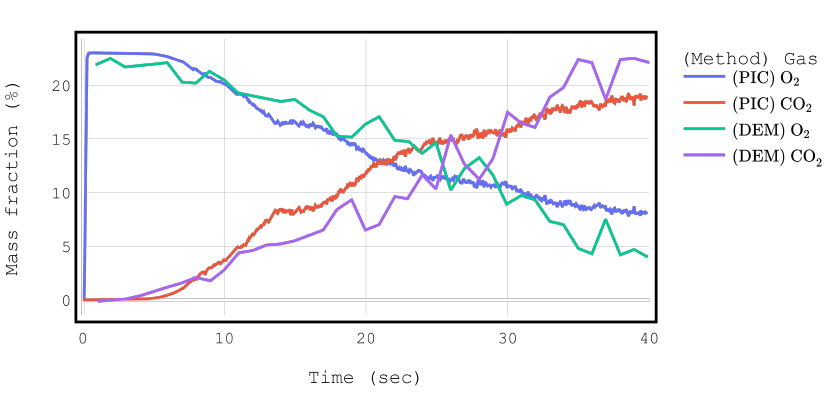}
\end{center}
\caption{Comparison of the $\mathrm{O}_2$ and $\mathrm{CO}_2$ mass fractions (\%) over 40s predicted by the PIC method (used in this work) and the discrete element method (see `DEM' in the figure) results obtained in~\citep{xie2021study}.}
\label{fig6}
\end{figure*} 
In \citep{XIE201972,xie2021study}, the discrete element method is used to compute the char combustion that couples heat transfer and complex chemical reactions for the boiler geometry configuration shown in Figure~\ref{fig5a}. he discrete element method provides highly accurate results as it simulates the individual particle behavior via detailed interactions and dynamics.  We thus validate our simulations with the results in \cite[Figure 2]{xie2021study} obtained by the discrete element method for the time evaluation of the gas mass fractions.
For the subsequent validation, we use the mean values of $X_1,X_2,\ldots,X_5$ from \cite[Table 1]{xie2021study}. 

Figure~\ref{fig6} presents the simulation results for the mass fractions of $\mathrm{CO}_2$ and $\mathrm{O}_2$ over $40$ seconds by the PIC method and the DEM results from \citep{xie2021study}. 
Both methods predict similar changes in the mass fraction of $\mathrm{CO}_2$ and $\mathrm{O}_2$. The $\mathrm{CO}_2$ mass fraction is inversely proportional to the $\mathrm{O}_2$ mass fraction. The Pearson correlation coefficients between the $\mathrm{CO}_2$ and $\mathrm{O}_2$ mass fractions during $1$--$40$ seconds are -0.9999 and -0.9975 from the results by PIC and DEM, respectively. These results indicate that both PIC and DEM demonstrate a strong inverse correlation between the mass fractions of $\rm{CO}_2$ and $\rm{O}_2$, implying that an increase in $\rm{CO}_2$ results in a decrease in $\rm{O}_2$.    
We compute the root-mean-squared error (RMSE) by measuring mass fractions at 5-second intervals over 40 seconds for both methods. The computed RMSE is 2.87 \%, indicating that the PIC model is accurate and well-validated compared to the discrete element method model.

\subsection{Time integrated quantity of interest}\label{sec:5.4}
For the global sensitivity analysis, we measure the total thermal energy of a mixed gas from zero to ten seconds as the QoI. 
We compute the thermal energy as 
\begin{align*} 
Q(\mathbf{X})=\int_{t=0}^{t=10}\dot{Q}(\mathbf{X},t)\mathrm{d}t=\int_{t=0}^{t=10}C_p(\mathbf{X},t) \times \dot{m}(\mathbf{X},t) \times  T_{\mathrm{avg}}(\mathbf{X},t)\mathrm{d}t,
\end{align*} 
where $T_{\mathrm{avg}}(\mathbf{X},t)$ is the average temperature across the outlet's cross-section (see $A_O$ in Figure~\ref{fig5a}) at time $t$ and inputs $\bX$, $\dot{m}(\bX,t)$ is the mass flow rate at the outlet and $C_{p}(\mathbf{X},t)$ is the specific heat capacity of the mixture, i.e., $C_p(\mathbf{X},t)=\sum_{i=1}^{i=5} C_{pi}(\mathbf{X},t)\times \textcolor{black}{M_i}$. Here, $C_{pi}$ is the specific heat capacity of each gas component in the mixture, and $M_i$ is the mole fraction of each gas component, i.e., $M_1=31.9988$ g/mol $(\rm{O}_2)$, $M_2=28.0134$ g/mol $(\rm{N}_2)$, $M_3=28.0104$ g/mol $(\rm{CO})$, $M_4=44.0098$ g/mol $(\rm{CO}_2)$, $M_5=18.0153$ g/mol $(\rm{H}_2\rm{O})$. 

\begin{figure*}[htbp]
\centering 
\begin{subfigure}{0.50\textwidth}
{\includegraphics[width=\linewidth]
{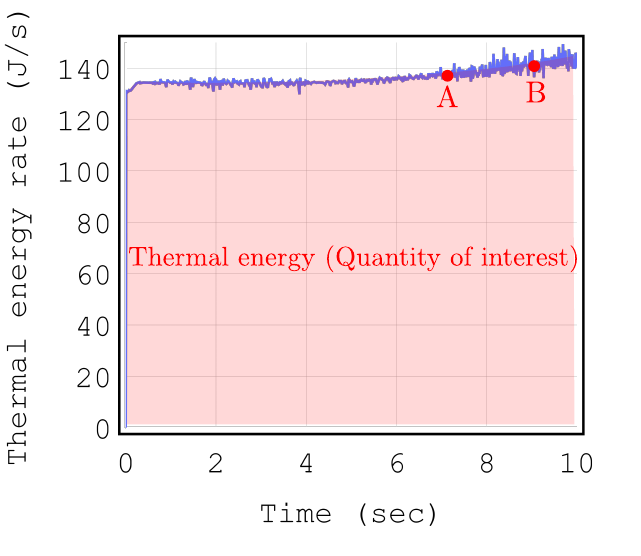}}\caption{Thermal energy}\label{fig7a}
\end{subfigure} 
\begin{subfigure}{0.43\textwidth}
{\includegraphics[width=\linewidth]
{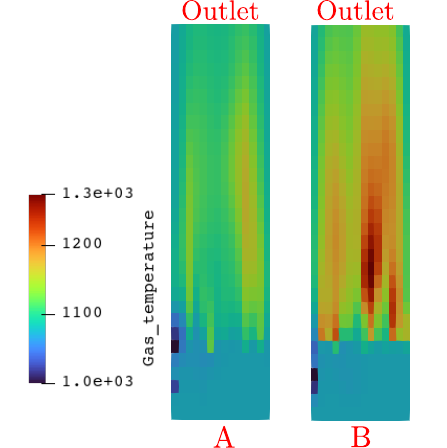}}\caption{Temperature}\label{fig7b}
\end{subfigure} 
\caption{(a) Thermal energy rate (J/s) measured at the outlet over 10 seconds when the mean values of the random inputs $(X_1,X_2,X_3,X_4,X_5)^{\intercal}$ are used. The measured total thermal energy over 10 seconds (highlighted in red color) is the QoI used for global sensitivity analysis. (b) Temperature distributions at 7 sec (A) and 9 sec (B) in units of Kelvin.}
\label{fig7}
\end{figure*}

Figure~\ref{fig7a} shows the time evolution of thermal energy rate measured at the outlet over 10 seconds. 
Figure~\ref{fig7b} presents the temperature distributions of the cross-section of the boiler at 7 and 9 seconds, as shown in panels A and B, respectively.   

\subsection{Surrogate model}\label{sec:5.5}
\begin{figure*}[htbp]
\centering 
\begin{subfigure}{0.47\textwidth}
{\includegraphics[width=\linewidth]
{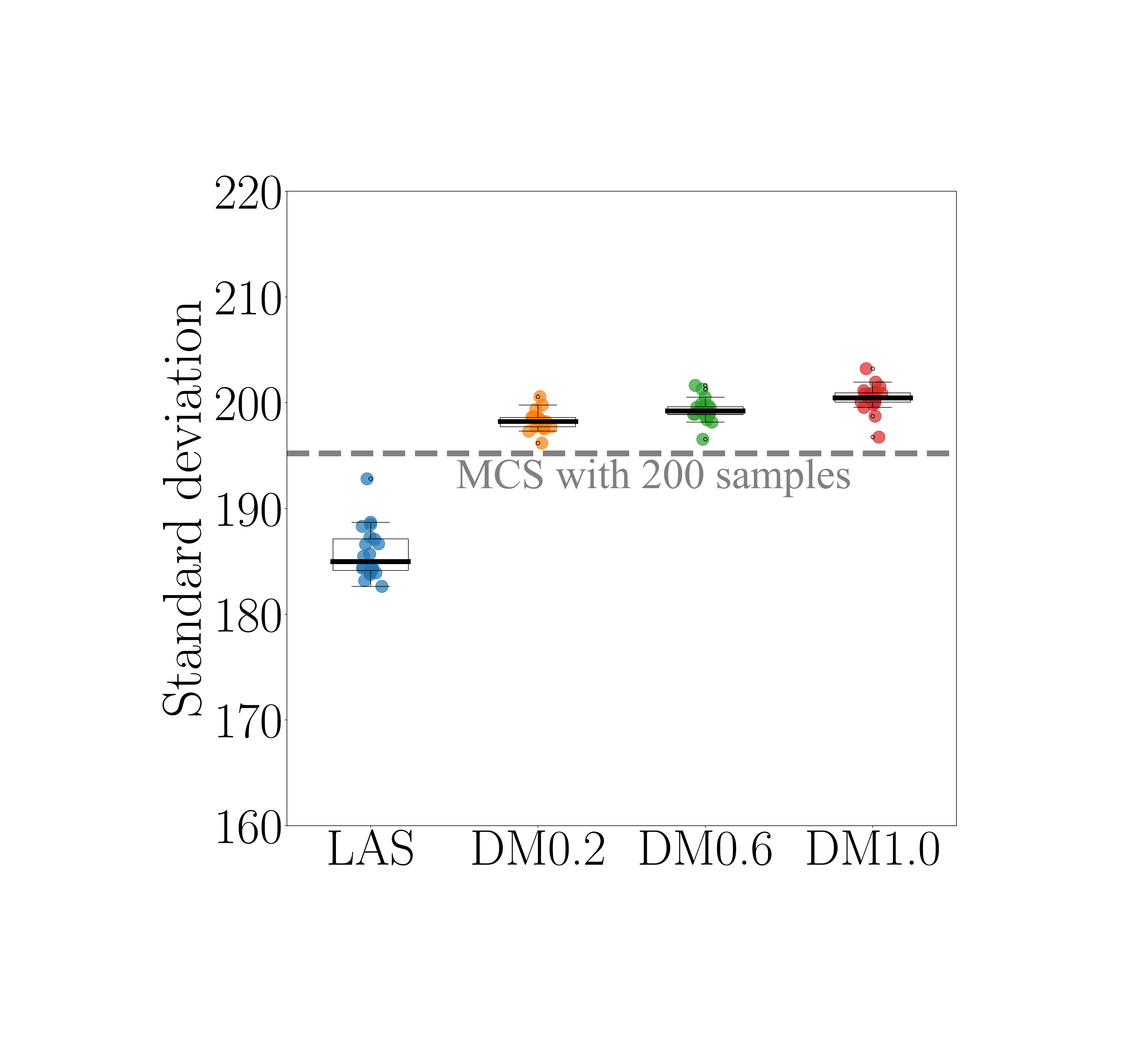}}\caption{Standard deviation}\label{fig8a}
\end{subfigure} 
\begin{subfigure}{0.47\textwidth}
{\includegraphics[width=\linewidth]
{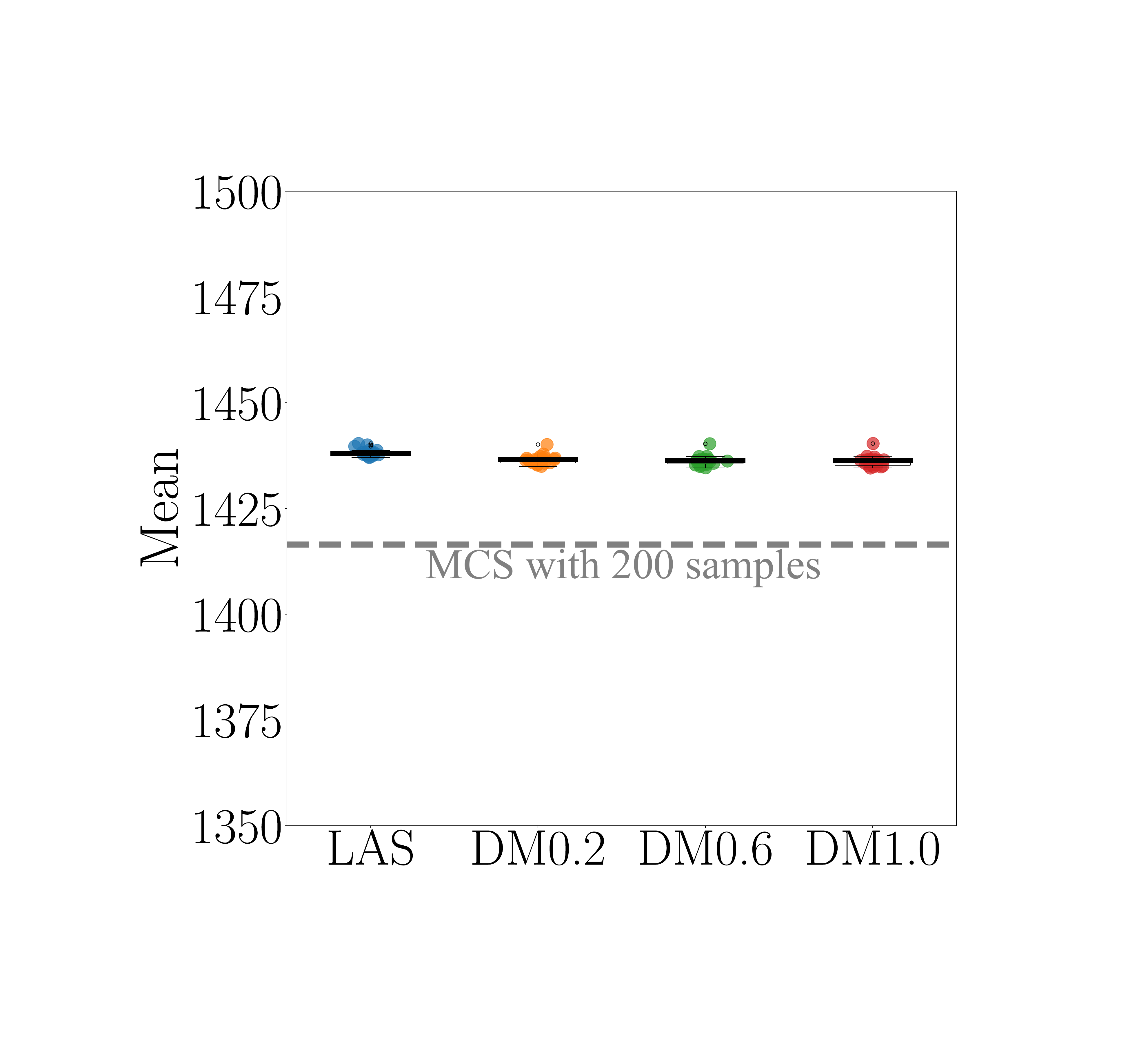}}\caption{Mean}\label{fig8b}
\end{subfigure} 
\caption{\color{black}{Boxplots of the standard deviation (a) and the mean (b) of the QoI, estimated by the bivariate eleventh-order PDD using the Lasso (LAS) and Lasso-based D-MORPH regressions with weight values $\lambda=0.2,0.6,1.0$ (DM0.2, DM0.6, DM1.0). Due to the budget limit (i.e., $M=195 < L=606$), we solve the underdetermined system using the Lasso and the proposed Lasso-based D-MORPH regression, repeating this process $20$ times. The unbiased Monte Carlo reference (MCS) with a sample size of 200 is shown as a gray-dotted line.}}
\label{fig8}
\end{figure*} 
\begin{figure*}[htbp]
\centering 
\begin{subfigure}{1.0\textwidth}
{\includegraphics[width=\linewidth]
{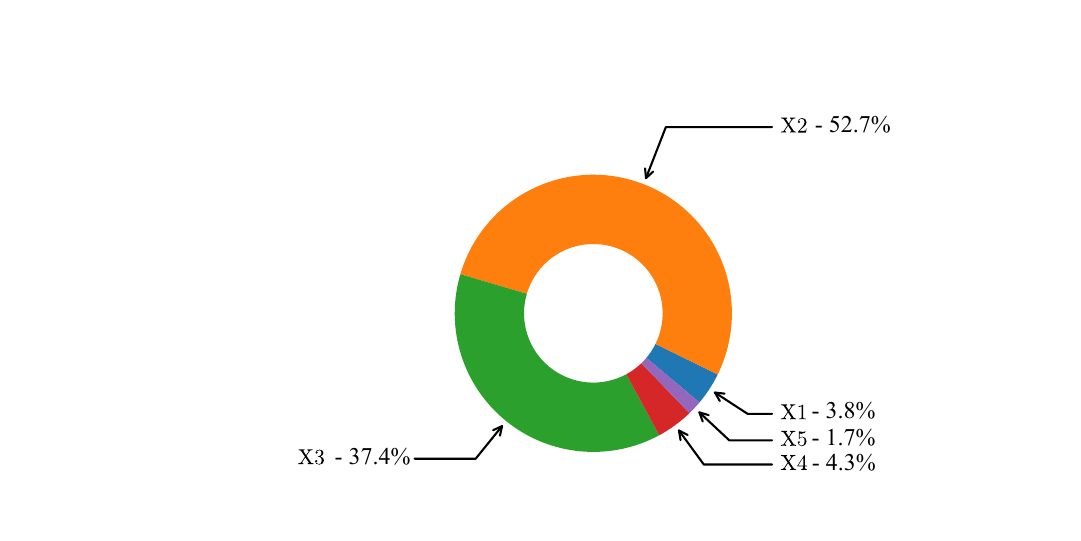}}\caption{Total effect sensitivities}\label{fig9a}
\end{subfigure} 
\begin{subfigure}{1.0\textwidth}
{\includegraphics[width=\linewidth]
{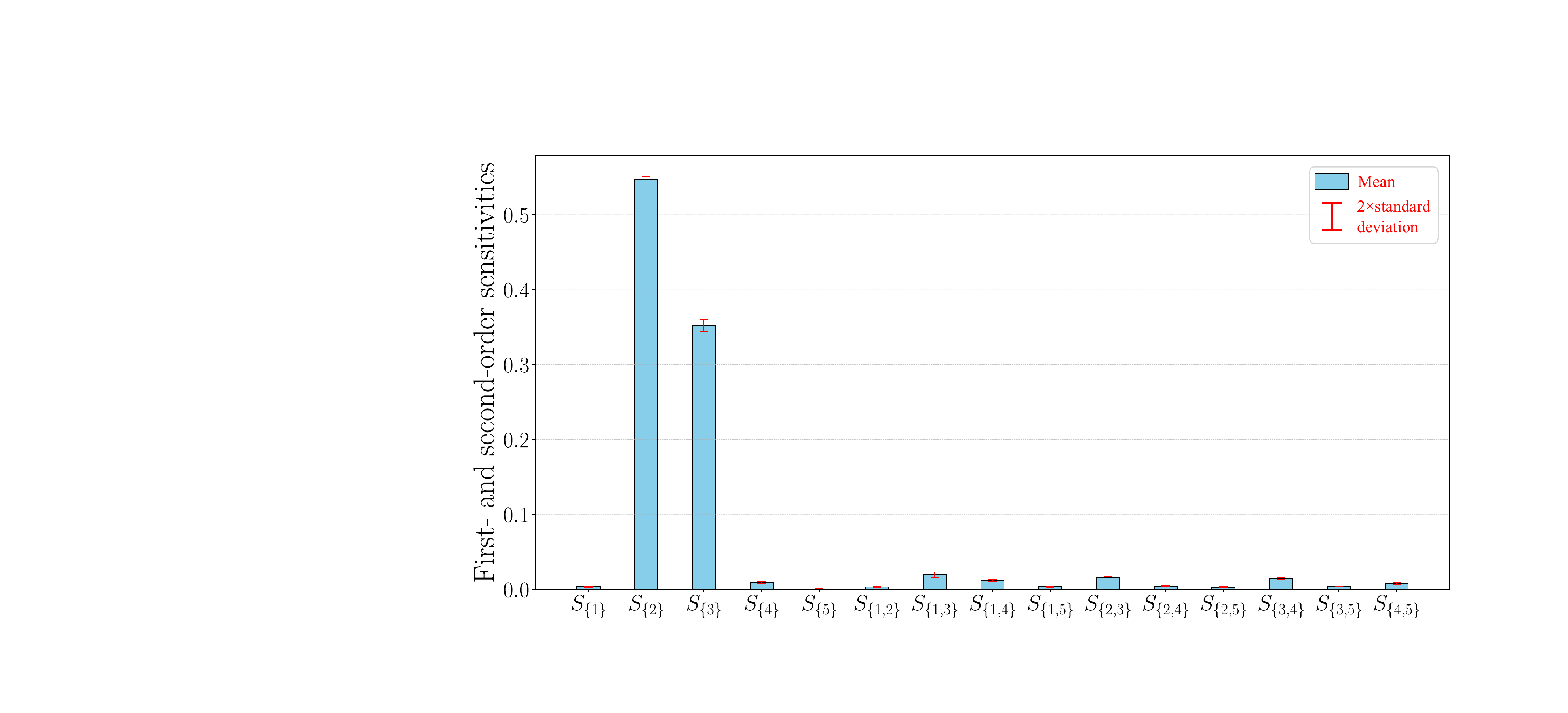}}\caption{First- and second-order sensitivities}\label{fig9b}
\end{subfigure} 
\caption{\textcolor{black}{(a) Total effect sensitivity indices using the PDD surrogate with the proposed D-MORPH regression and weight $\lambda=0.2$. This result shows that $X_2$ (air inflow) and $X_3$ (diameter of char particles) are the most influential to the QoI, both each close to 52.7\% and 37.4\%, whereas $X_5$ (char mass inflow) is the least influential. The other inputs $X_1$ (freeboard height) and $X_4$ (diameter of the glass bead particles) have 4\% impact on the QoI and (b) estimates of first- and second-order sensitivity indices}}
\label{fig9}
\end{figure*}  
We consider 
a computational budget of 1,545 CPU hours. Within this limit, we obtain \textcolor{black}{195 training samples for computing the PDD surrogate (i.e., $B$=1,545 hours > 1,540.5 hours = 195 samples $\times$ 7.9 hours/sample}). 
Given the PDD parameters $S=2$ and \textcolor{black}{$m=11$}, the surrogate model includes $606$ basis functions; we determine the corresponding $606$ unknown expansion coefficients via a regression method. With $M=195$ training samples representing 32\% of the $L=606$ unknown expansion coefficients (leading to an underdetermined linear system in \eqref{linear}), we use the proposed D-MORPH regression and a standard Lasso regression.

Figure~\ref{fig8a} presents boxplots of the standard deviation estimates obtained from a PDD surrogate computed by the proposed D-MORPH regression and a standard Lasso regression over 20 experiments. 
This figure also includes an unbiased reference obtained by Monte Carlo simulation with 200 samples, shown as in a gray dotted line.
For the standard deviation, the proposed D-MORPH results \textcolor{black}{(198.26, 199.30, 200.42)} of the three weight cases $\lambda$=0.2, 0.6, 1.0, averaged over 20 experiments, are between $2-3$x closer to the reference \textcolor{black}{(195.24)} by Monte Carlo simulation than the Lasso result (185.82). The standard deviation results of the D-MORPH regression with weights $\lambda$=0.2, 0.6, 1.0 are also close to each other.    
\textcolor{black}{For the statistical mean, the D-MORPH regression estimates \textcolor{black}{(1436.56, 1436.24, 1436.13)} for the three weight cases $\lambda$=0.2, 0.6, 1.0, averaged over 20 experiments, show slightly closer alignment to the reference ($1414.64$) than the LAS estimate ($1438.20$).} 

Overall, the proposed D-MORPH regression performs better than the original Lasso regression when estimating the standard deviation. We select the weight \textcolor{black}{$\lambda=0.2$} (DM0.2) to implement global sensitivity analysis, \textcolor{black}{as it produces standard deviation estimates with the least bias}. However, we note that the number of Monte Carlo samples \textcolor{black}{(200)} is not enough for accurate benchmark estimates.

\subsection{Global sensitivity analysis}  \label{sec:5.6}
We perform a global sensitivity analysis for the char combustion model via its surrogate described above. Figure~\ref{fig9} presents the global sensitivity analysis results obtained by the PDD surrogate with the proposed D-MORPH regression.

\textcolor{black}{In Figure~\ref{fig9a}, the pie chart presents the total effect sensitivities, as a portion of 100\%, for inputs $X_1$ through $X_5$. The chart demonstrates that $X_2$ (air inflow) has the most influence (52.7\%) on the QoI, while $X_5$ (char mass inflow) has the least influence (1.7\%).
The $X_3$ (char diameter) also has a significant impact ($37.4\%$) on the QoI, while each of $X_1$ (freeboard height) and $X_4$ (diameter of glass bead) contributes a relatively small proportion (4\%) to the total effect sensitivities.} 

\textcolor{black}{In Figure~\ref{fig9b}, the bar chart reveals that the first-order sensitivities for inputs $X_2$ and $X_3$ are the most influential, each being larger than 54\% and 35\%, respectively.}  
\textcolor{black}{In contrast, the second-order sensitivities for all inputs are small. Although the first-order sensitivities $S_{\{1\}} (0.36\%)$ and $S_{\{4\}} (1\%)$ are almost negligible, their second-order sensitivities indicate that $S_{\{1\}}$ and $S_{\{4\}}$ have a higher impact, resulting in $4\%$ each in total sensitivities, as shown in Figure~\ref{fig9a}.   Therefore, the Sobol sensitivity method provides a more detailed picture of the sensitivity contributions than local sensitivity and Morris screening~\citep{morris1991factorial}, both of which only capture first-order indices.}

This is important information for system designers and operators. When one optimizes the boiler's thermal energy, the results from this global sensitivity analysis allow one to prioritize the sensitive inputs, i.e., \textcolor{black}{air inflow and char diameter over glass bead diameter, char mass inflow, and freeboard height.} 

\section{Conclusions} \label{sec:6} 
%
We developed a novel D-MORPH regression method to train a PDD surrogate of the quantity of interest in the underdetermined case, which is common in applications where only limited data is available. The new regression method integrates a sparse Lasso solution into the cost function. The difference between the D-MORPH solution and the Lasso solution is minimized iteratively by a D-MORPH regression process. 
We evaluated the proposed method for global sensitivity analysis in two mathematical functions. A nonlinear function of five random inputs showed that the iterative algorithm improved and converged the D-MORPH solution. In the second mathematical function, we considered fifteen (high-dimensional) random inputs. The results clearly showed that the proposed D-MORPH regression is 3$\times$ times more accurate than the standard Lasso regression, and the proposed D-MORPH-based estimates for sensitivities converge as the training sample number increase. 
We finally demonstrated the versatility and scalability of the PDD surrogate with the proposed D-MORPH regression for global sensitivity analysis of a lab-scale boiler for char combustion. We only used 151 training samples, which are 15\% of the data required when using a standard regression method. Consequently, we found four inputs that influence the QoI (thermal energy), showing that the other input is relatively less influential. 

This work assumes that input random variables are independent. In practice, there are cases where input variables are dependent or correlated. A potential approach to extend the proposed method to these cases would require the use and development of other surrogates, such as generalized polynomial chaos expansion~\citep{lee2020practical,jakeman2019polynomial} or dimensionally decomposed generalized polynomial chaos expansion~\citep{lee2023bi,lee2023high}. This surrogate is more accurate and efficient in handling dependent inputs as the surrogate uses orthonormal polynomial bases for dependent random inputs. Performing a global sensitivity analysis with dependent inputs would then also require a covariance decomposition method instead of the variance-based method.

\section*{Acknowledgments} \label{sec:7} 
This research was in part financially supported by the Korea Institute for Advancement of Technology (KIAT) through the International Cooperative R\&D program (No. P0019804, Digital twin based intelligent unmanned facility inspection solutions). We thank Dr. Kyung Man Kim from Virnect for helpful discussions about the char combustion application. 

\bibliography{main}
\bibliographystyle{abbrv}
\newpage

\newpage
\begin{appendices}
\section{Sensitivity results for random inputs $X_7$--$X_{15}$ in example~2\label{sec:appx1}}
\begin{figure*}[ht]
\centering 
\begin{subfigure}{0.495\textwidth}{\includegraphics[width=\linewidth, trim={8cm 2cm 8cm 6cm}, clip]{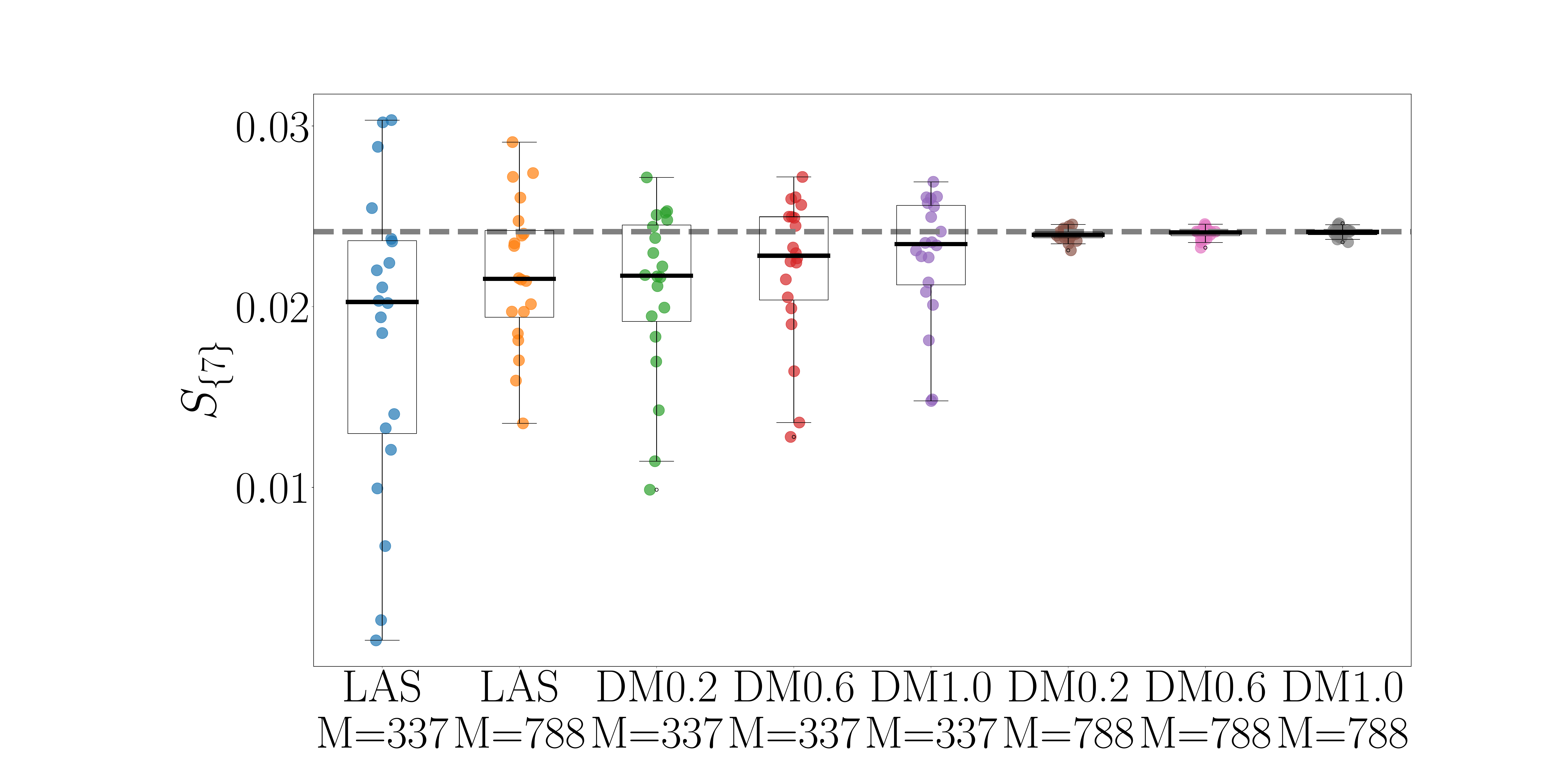}}\caption{$S_{\{7\}}$}\label{fig10:(a)}
\end{subfigure} 
\begin{subfigure}{0.495\textwidth}{\includegraphics[width=\linewidth, trim={8cm 2cm 8cm 6cm}, clip]{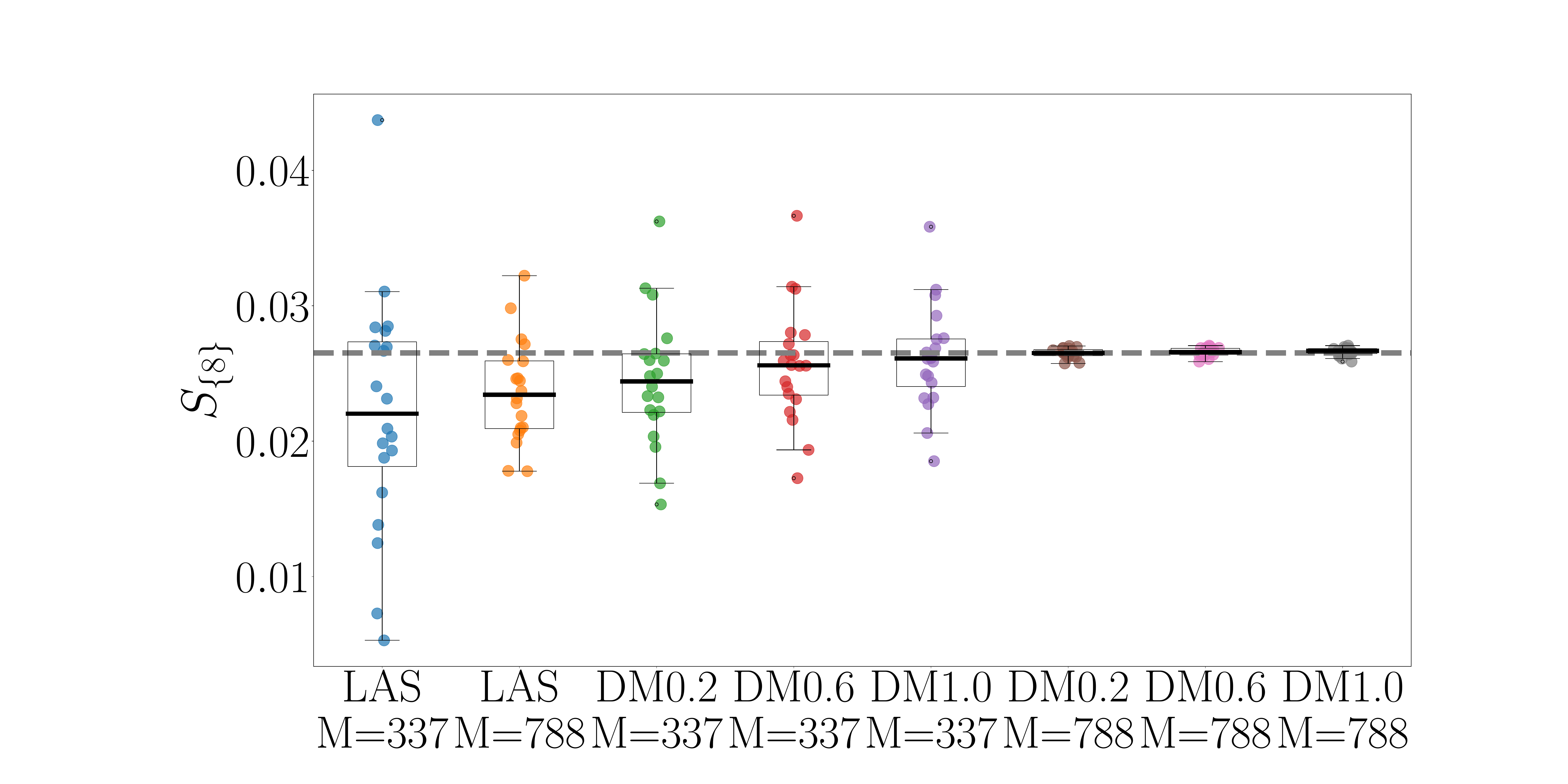}}\caption{$S_{\{8\}}$}\label{fig10:(b)}
\end{subfigure} 
\begin{subfigure}{0.495\textwidth}{\includegraphics[width=\linewidth, trim={8cm 2cm 8cm 6cm}, clip]{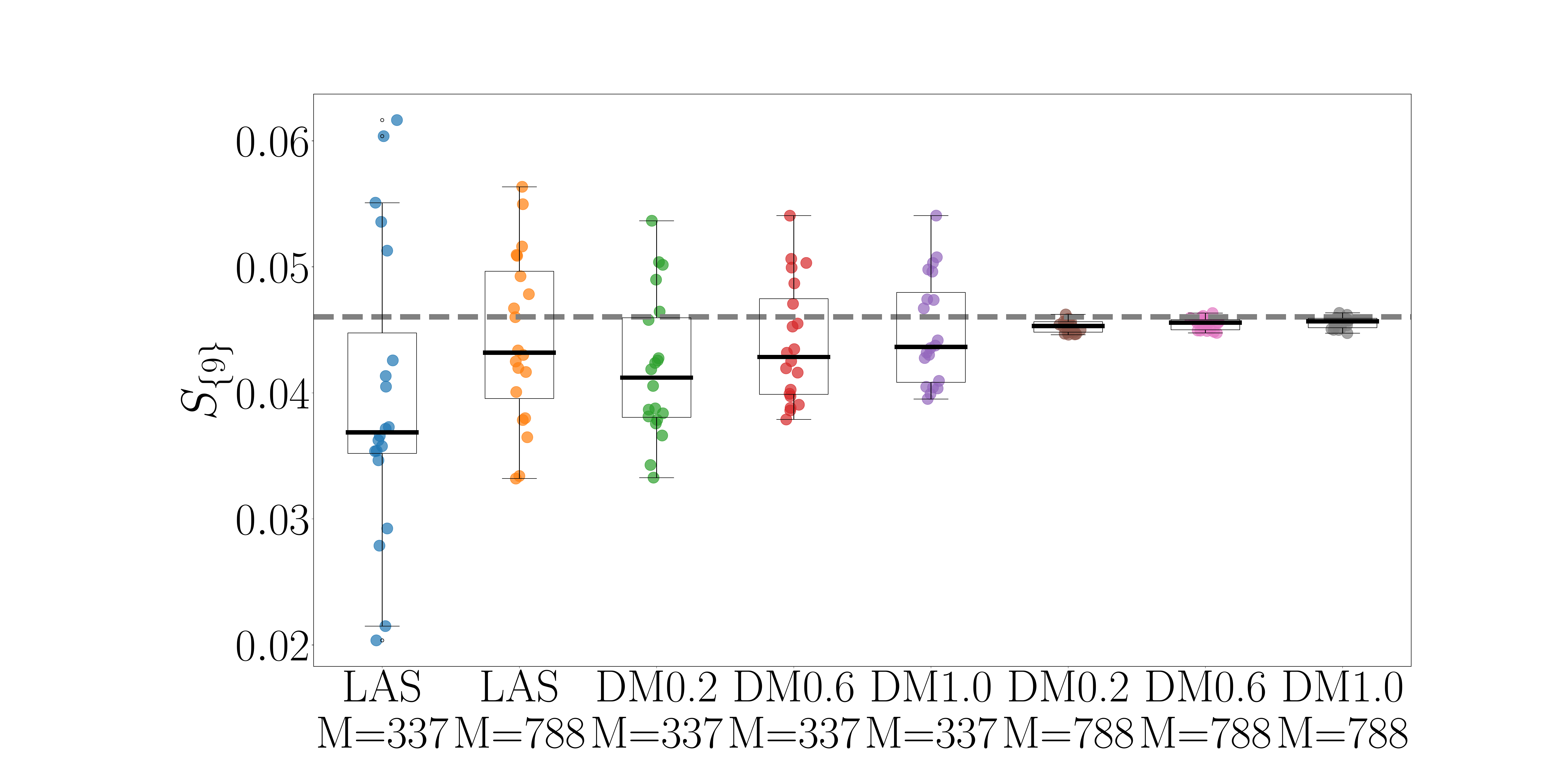}}\caption{$S_{\{9\}}$}\label{fig10:(c)}
\end{subfigure} 
\begin{subfigure}{0.495\textwidth}{\includegraphics[width=\linewidth, trim={8cm 2cm 8cm 6cm}, clip]{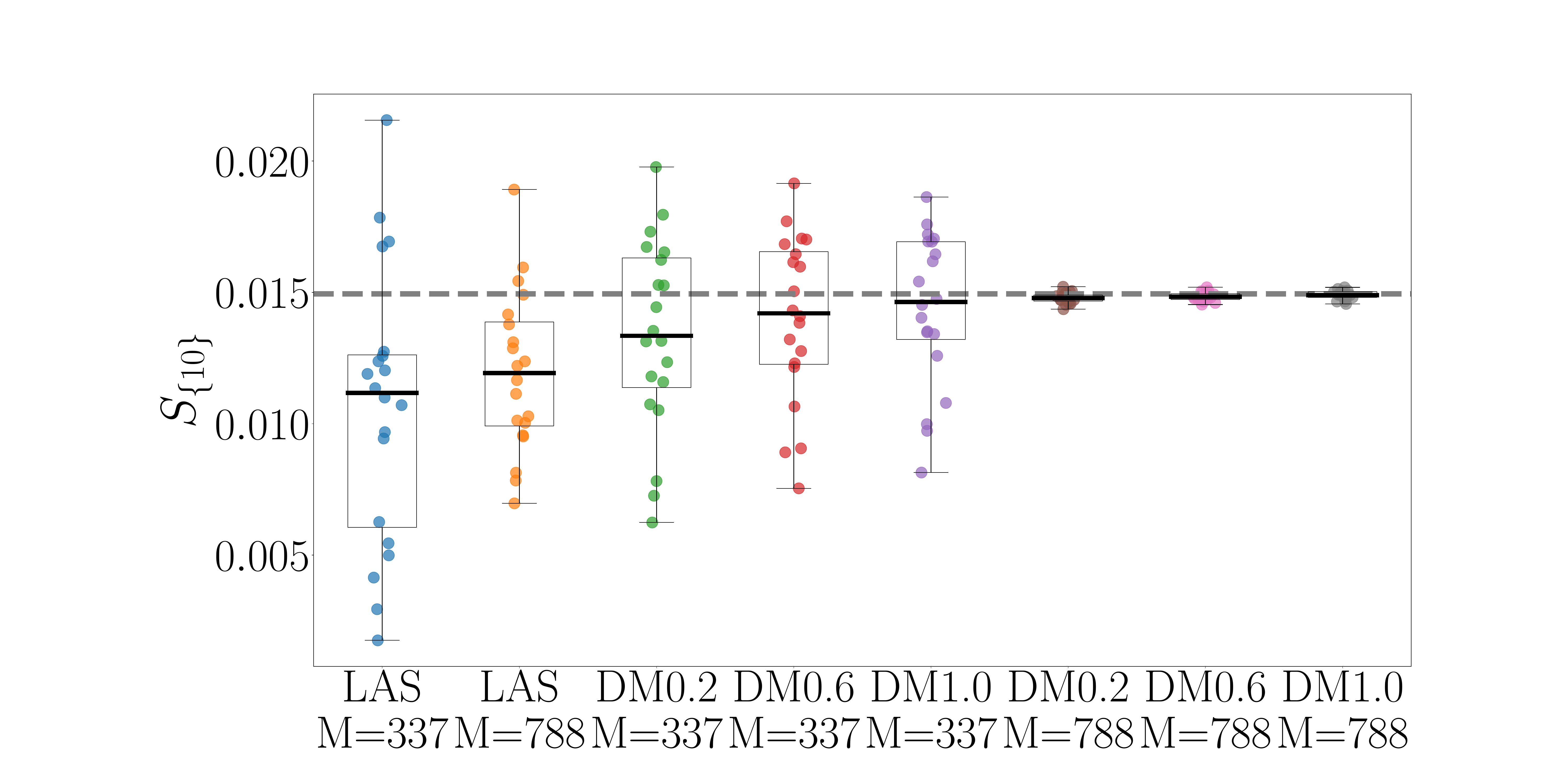}}\caption{$S_{\{10\}}$}\label{fig10:(d)}
\end{subfigure} 
\begin{subfigure}{0.495\textwidth}{\includegraphics[width=\linewidth, trim={8cm 2cm 8cm 6cm}, clip]{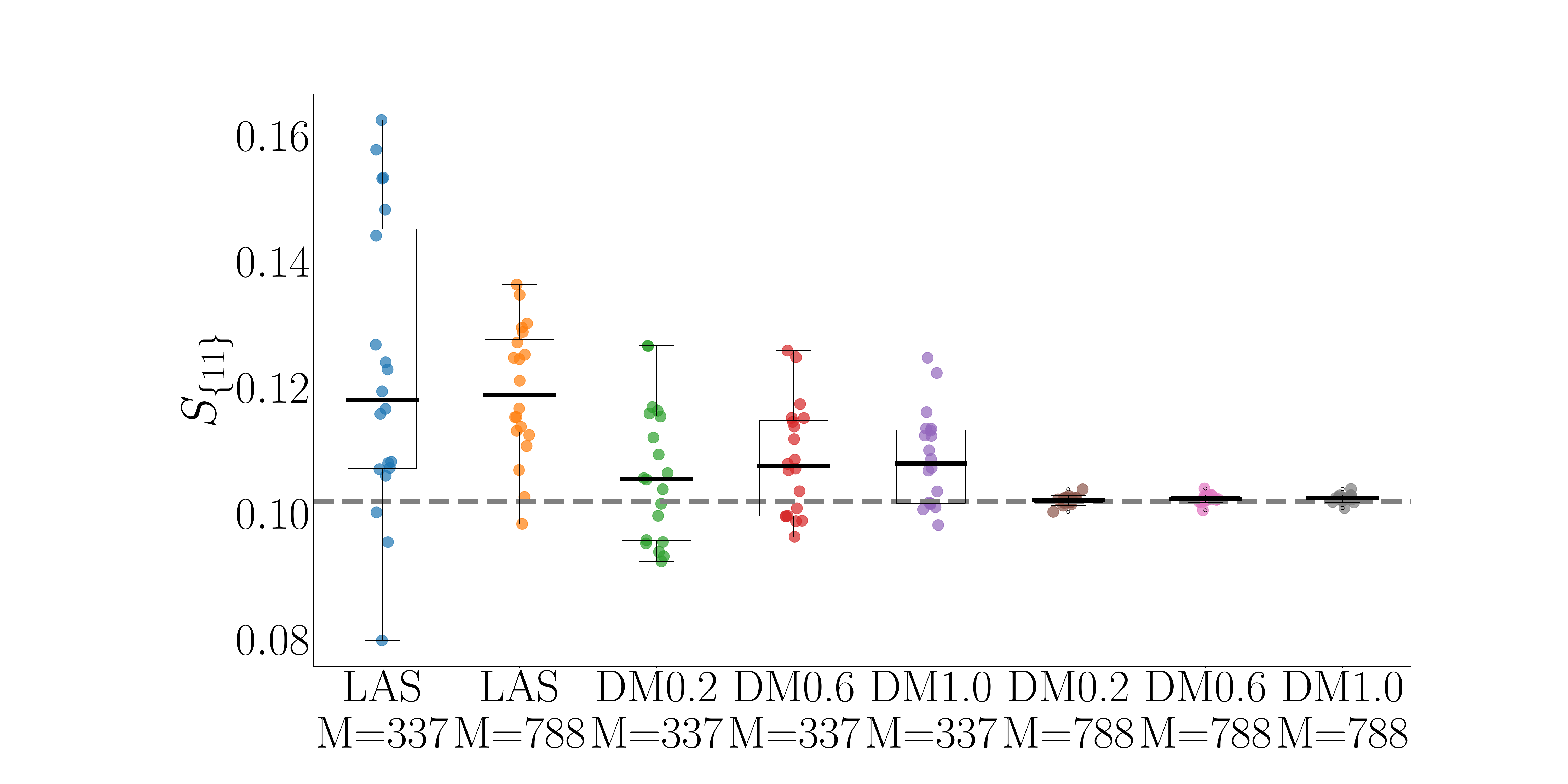}}\caption{$S_{\{11\}}$}\label{fig10:(e)}
\end{subfigure}
\begin{subfigure}{0.495\textwidth}{\includegraphics[width=\linewidth, trim={8cm 2cm 8cm 6cm}, clip]{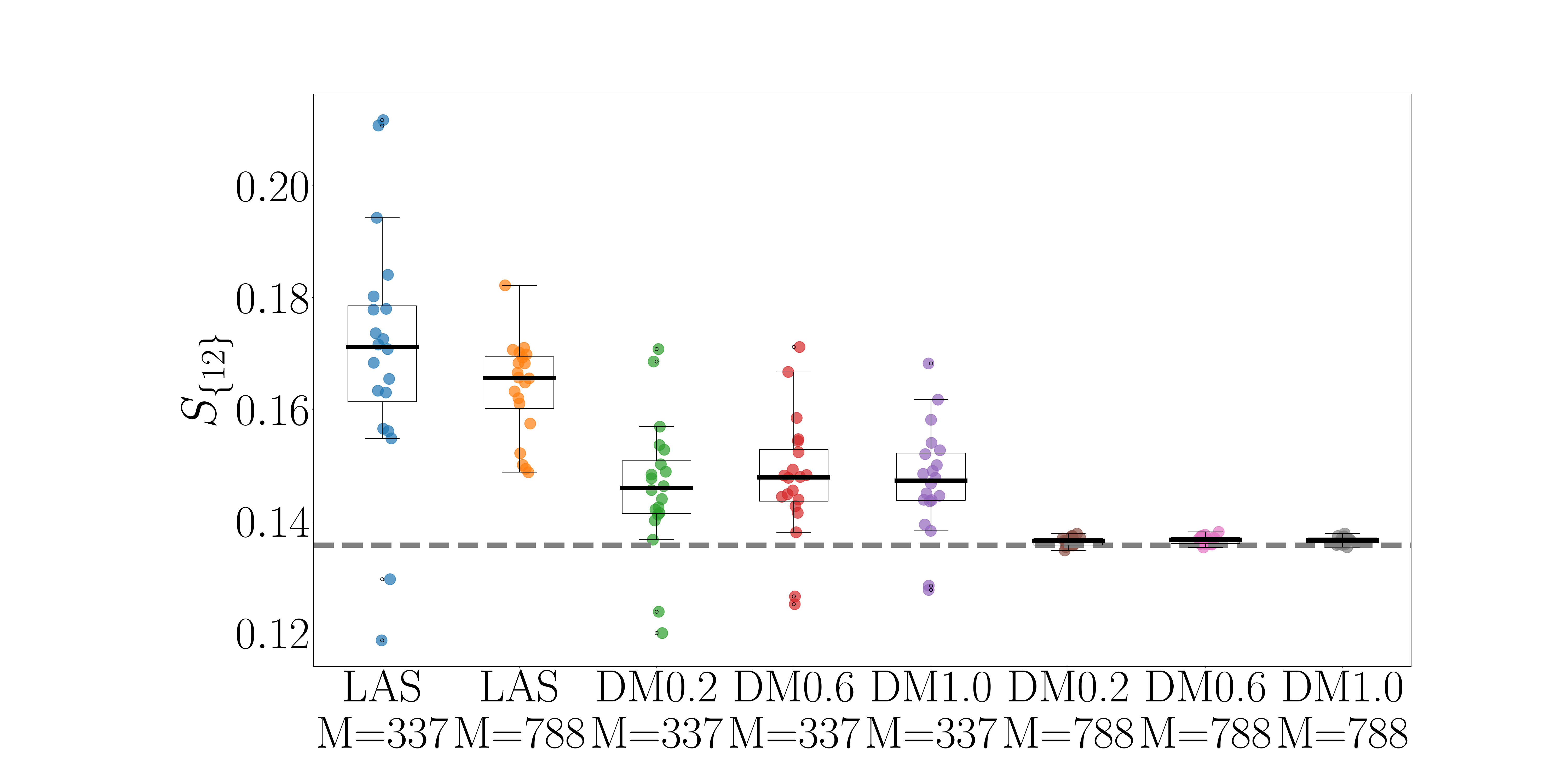}}\caption{$S_{\{12\}}$}\label{fig10:(f)}
\end{subfigure}
\caption{Boxplots of first-order sensitivity indices $S_{\{i\}}$, $i=7,8,9,10,11,12$, of the random output $y(\bX)$, in (a), (b), (c), (d), (e), (f), respectively, estimated by the bivariate fifth-order PDD using the Lasso (LAS) and Lasso-based D-MOPRH regressions with weight values $\lambda=0.2$, $0.6$, $1.0$ (DM0.2, DM0.6, DM1.0). Two underdetermined systems are considered: $M=337$ and $M=788$, which correspond to 30\% and 70\% of the number ($L=1,126$) of expansion coefficients. Each regression is repeated 20 times.  The exact solution is shown as a gray-dotted line.}
\label{fig10}
\end{figure*} 
\begin{figure*}[ht]
\centering 
\begin{subfigure}{0.495\textwidth}{\includegraphics[width=\linewidth, trim={8cm 2cm 8cm 6cm}, clip]{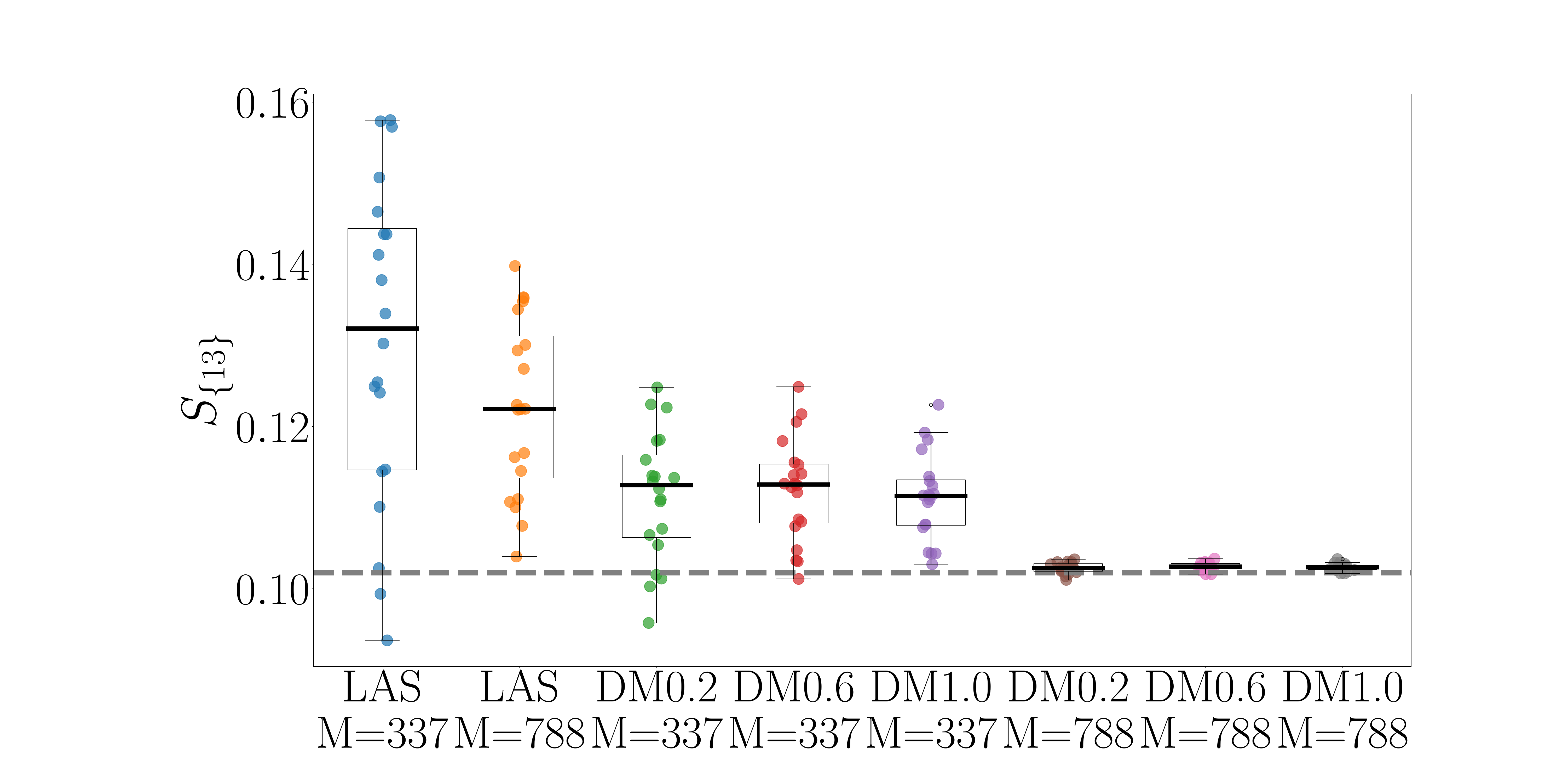}}\caption{$S_{\{13\}}$}\label{fig11:(a)}
\end{subfigure} 
\begin{subfigure}{0.495\textwidth}{\includegraphics[width=\linewidth, trim={8cm 2cm 8cm 6cm}, clip]{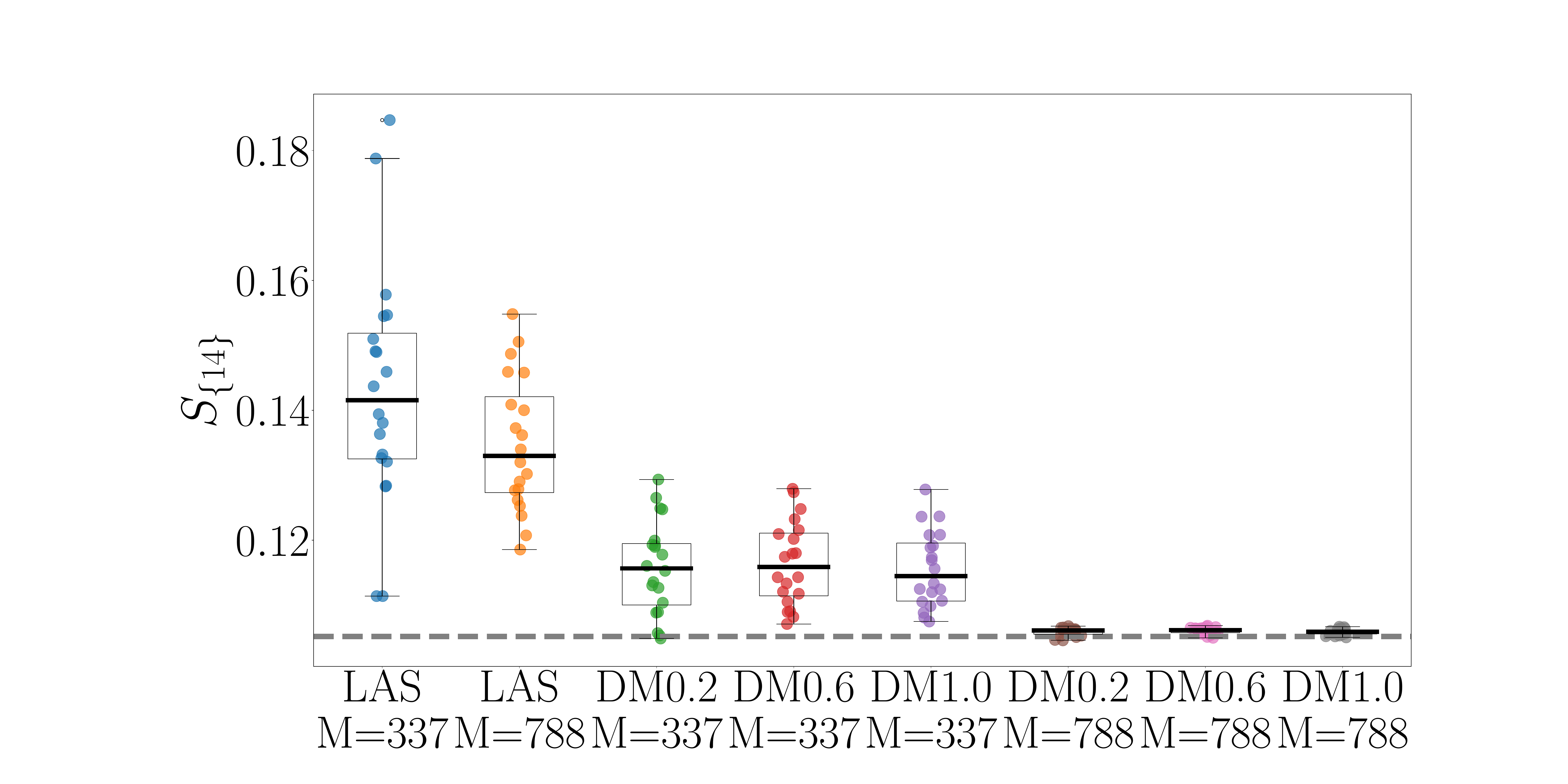}}\caption{$S_{\{14\}}$}\label{fig11:(b)}
\end{subfigure} 
\begin{subfigure}{0.495\textwidth}{\includegraphics[width=\linewidth, trim={8cm 2cm 8cm 6cm}, clip]{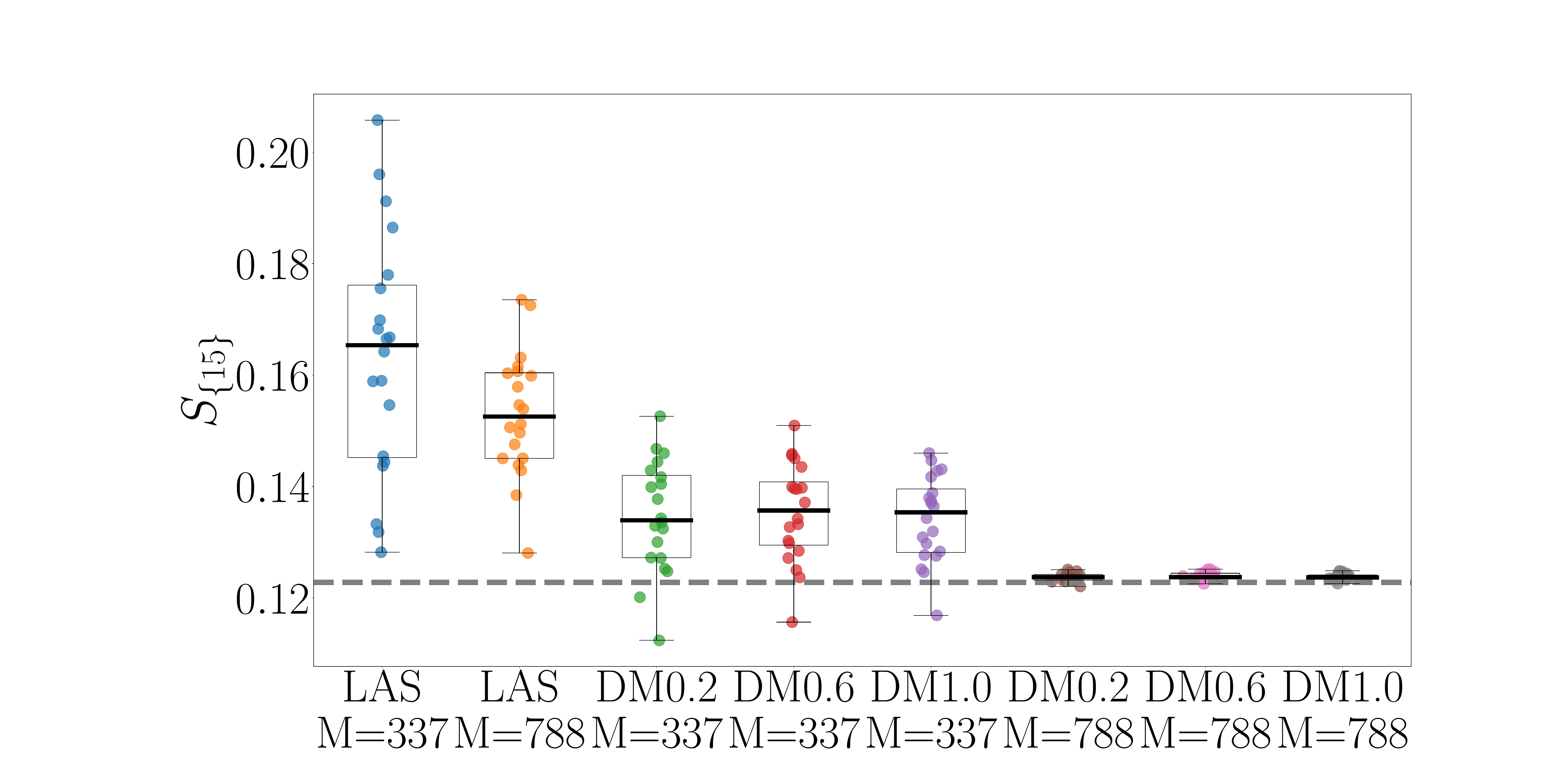}}\caption{$S_{\{15\}}$}\label{fig11:(c)}
\end{subfigure} 
\caption{Boxplots of first-order sensitivity indices $S_{\{i\}}$, $i=13,14,15$, of the random output $y(\bX)$, in (a), (b), (c), respectively, estimated by the bivariate fifth-order PDD using the Lasso (LAS) and Lasso-based D-MOPRH regressions with weight values $\lambda=0.2$, $0.6$, $1.0$ (DM0.2, DM0.6, DM1.0). Two underdetermined systems are considered: $M=337$ and $M=788$, which correspond to 30\% and 70\% of the number ($L=1,126$) of expansion coefficients. Each regression is repeated 20 times.  The exact solution is shown as a gray-dotted line.}
\label{fig11}
\end{figure*} 
\end{appendices} 
\end{document}